\documentclass[draft]{amsart}
\usepackage{amsmath,amssymb,amsthm,eucal,verbatim}

\newtheorem{theorem}{Theorem}
\newtheorem*{theorem*}{Theorem}
\newtheorem{conjecture}{Conjecture}
\newtheorem*{conjecture*}{Conjecture}
\newtheorem{lemma}{Lemma}

\newtheorem*{corollary*}{Corollary}

\newtheorem*{HypoB}{Hypothesis $BV(\theta)$}
\newtheorem*{TheoB}{Theorem (Bombieri-Vinogradov)}
\newtheorem*{TheoB2}%
{Theorem (Bombieri-Vinogradov for $\varpi*\varpi$)}
\newtheorem*{HypoB2}{Hypothesis $BV_2(\theta)$}

\numberwithin{equation}{section}
\newcommand{\cE}{{\mathcal E}}
\newcommand{\cH}{{\mathcal H}}
\newcommand{\cL}{{\mathcal L}}
\newcommand{\cM}{{\mathcal M}}
\newcommand{\cS}{{\mathcal S}}
\newcommand{\bZ}{{\mathbb Z}}
\newcommand{\bb}{{\bf b}}
\newcommand{\bM}{{\bf M}}
\newcommand{\ctS}{\tilde{\cS}}
\newcommand{\tb}{\tilde{b}}

\newcommand{\fS}{\mathfrak S}
\newcommand{\asymptotic}{\sim}
\newcommand{\binomial}{\binom}
\newcommand{\union}{\cup}

\newcommand{\SEP}{\beta(\cH)\fS(\cH)/\log R}
\newcommand{\SE}{\left\{1+O(\SEP)\right\}}

\def\sumprime_#1{\setbox0=\hbox{$\scriptstyle{#1}$}
\setbox2=\hbox{$\displaystyle{\sum}$}
\setbox4=\hbox{${}'
\mathsurround=0pt$}
\dimen0=.5\wd0 \advance\dimen0 by-.5\wd2
\ifdim\dimen0>0pt
\ifdim\dimen0>\wd4 \kern\wd4 \else\kern\dimen0\fi\fi
\mathop{{\sum}'}_{\kern-\wd4 #1}}

\def\sumflat_#1{\setbox0=\hbox{$\scriptstyle{#1}$}
\setbox2=\hbox{$\displaystyle{\sum}$}
\setbox4=\hbox{${}\flat
\mathsurround=0pt$}
\dimen0=.5\wd0 \advance\dimen0 by-.5\wd2
\ifdim\dimen0>0pt
\ifdim\dimen0>\wd4 \kern\wd4 \else\kern\dimen0\fi\fi
\mathop{{\sum}^\flat}_{\kern-\wd4 #1}}

\begin{document}
\title{Small gaps between  primes or almost primes}

\author{D. A. Goldston}
\address{Department of Mathematics, San Jose
State University, San Jose, CA 95192, USA}
\email{goldston@math.sjsu.edu}
\author
{S.W. Graham}
\address
{Department of Mathematics, Central Michigan University,
Mt. Pleasant, MI 48859, USA}
\email{graha1sw@cmich.edu}
\author{J. Pintz}
\address{R\'enyi Mathematical Institute of the Hungarian Academy
of Sciences, H-1364 Budapest, P.O.B. 127, Hungary} 
\email{pintz@renyi.hu}
\author{C. Y. Y{\i}ld{\i}r{\i}m}
\address{Department of Mathematics, Bo\~{g}azi\c{c}i University,
Istanbul 34342 \& \newline 
Feza G\"{u}rsey Enstit\"{u}s\"{u}, \c{C}engelk\"{o}y, Istanbul,
P.K. 6, 81220, Turkey}
\thanks{The first author was supported by NSF grant DMS-0300563,
the NSF Focused Research Group grant 0244660, and the 
American Institute of Mathematics;
the second author by a sabbatical leave from 
Central Michigan University; 
the third author by OTKA grants No. T38396, T43623, T49693 and 
the Balaton program; 
the fourth author by T\"UB\.{I}TAK}

\date{June 2, 2005}

\begin{abstract}
Let $p_n$ denote the $n^{{\rm th}}$ prime. Goldston, Pintz, and Yildirim
recently proved that 
\begin{equation*}
    \liminf_{n\to \infty} \frac{(p_{n+1}-p_n)}{\log p_n} =0.
\end{equation*}
We give an alternative proof of this result. We also prove some
corresponding results for numbers with two prime factors. 
Let $q_n$ denote the $n^{{\rm th}}$ number that is a product of 
exactly two distinct primes. We prove that
\begin{equation*}
    \liminf_{n\to \infty} (q_{n+1}-q_n) \le 26.
\end{equation*}
If an appropriate generalization of the Elliott-Halberstam Conjecture
is true, then the above bound can be improved to $6$.

\end{abstract}

\maketitle

\section{Introduction}

In 1849, A. de Polignac (\cite{D}, p. 424) 
conjectured that every even number is
the difference of two primes in infinitely many ways. 
More generally, we can let 
$\cH=\{h_1, h_2, \ldots, h_k\}$ be a set of $k$ distinct
integers. 
A major open question in number theory is to show
that there are infinitely many positive integers
$n$ such that $n+h_1, n+h_2, \ldots, n+h_k$ are all 
prime, provided that $\cH$ meets an obvious necessary
condition that we call {\it admissibility}.
For each prime $p$, let $\nu_p(\cH)$ be the number of 
distinct residue classes mod $p$ in $\cH$. We say that the set 
$\cH$ is {\it admissible} if $\nu_p(\cH)<p$ for all $p$.

Using heuristics from the circle method,
Hardy and Littlewood \cite{HL3} realized the significance
of the singular series $\fS(\cH)$, defined as  
\begin{equation} \label{E:fSDefinition}
    \fS(\cH)=\prod_{p} \left(1-\frac{\nu_p(\cH)}{p} \right)
                  \left(1-\frac{1}{p} \right)^{-k} 
\end{equation}
for this problem. They made a conjecture about 
the asymptotic distribution of the numbers $n$ 
for which $n+h_1, \ldots , n+h_k$ are all prime,
which we state here in the following form.

\begin{conjecture}
    Let $\varpi(n)$ denote function
    \begin{equation} \label{E:ThetaDef}
	\varpi(n)=
	\begin{cases} 
	    \log n & \text{ if $n$ is prime,}\\
	    0      & \text{ otherwise.}
	\end{cases}
    \end{equation}
    As $N$ tends to infinity,
\begin{equation} \label{E:HLConjecture}
    \sum_{n\le N} \varpi(n+h_1) \varpi(n+h_2) \ldots \varpi(n+h_k)
    =N (\fS(\cH) +o(1)).
\end{equation}
\end{conjecture}

From the definition of $\fS(\cH)$, we see that $\fS(\cH)\ne 0$
if and only if $\nu_p(\cH)<p$ for all primes $p$; i.e.,
if and only if $\cH$ is admissible.  

The set $\cH=\{0,2\}$ is admissible, so the Hardy-Littlewood 
conjecture implies that 
\begin{equation*}
    \liminf_{n\to \infty}( p_{n+1} - p_n ) =2,
\end{equation*}
where $p_n$ denotes the $n^{\rm th}$ prime. 
In an unpublished paper in the {\it Partitio Numerorum}
series, Hardy and Littlewood \cite{HL7} proved 
that if the Generalized Riemann Hypothesis is true,
then
\begin{equation*}
    \liminf_{n\to \infty}
    \left( 
       \frac{ p_{n+1}-p_n }{\log p_n} 
      \right) \le \frac23.
\end{equation*}

In 1940, Erd\H{o}s \cite{Erdos} used Brun's sieve to give the first
unconditional proof of the inequality
\begin{equation*}
    \liminf_{n\to \infty}
    \left( 
       \frac{ p_{n+1}-p_n }{\log p_n} 
      \right) <1.
\end{equation*}
In 1965, Bombieri and Davenport \cite{BD} proved unconditionally that
\begin{equation} \label{E:BDBound}
    \liminf_{n\to \infty}
    \left( 
       \frac{ p_{n+1}-p_n }{\log p_n} 
      \right) \le 0.4665 \ldots.
\end{equation}
This result was one of the first applications of what is now known as 
the ``Bombieri-Vinogradov Theorem,'' which we state as follows.

\begin{TheoB} When $(a,q)=1$, let $E(x;q,a)$ be defined by the relation
\begin{equation} \label{E:EDef}
\sum_{\substack{x< n \le 2x \\ n\equiv a\pmod q}} 
   \varpi(n) = 
   \frac{x}{\phi(q)} + E(x;q,a).
\end{equation}
Furthermore, let 
\begin{equation} \label{E:E*Def}
    E(x,q) = \max_{a; (a,q)=1} |E(x,q,a)|,\quad
    E^*(N,q) = \max_{x\le N} E(x,q).
\end{equation}
If $A> 0$, then there exists $B>0$ such that if  
$Q\le N^{1/2} \log^{-B} N$, then
\begin{equation} \label{E:TheoBConclusion}
\sum_{q\le Q} E^*(N,q)
\ll_{A} N(\log N)^{-A}.
\end{equation}
\end{TheoB}

This result was proved by Bombieri in 1965 (\cite{Bo}). At about the 
same time, A. I. Vinogradov (\cite{Vi}) gave an independent proof of a 
slightly weaker result. There are numerous proofs of this result 
available in the literature; see, for example, \cite{Da} and  
\cite{Va}. We remark that in the usual definition of $E(x;q,a)$, one
takes the sum in \eqref{E:EDef} to be over $n\le x$. However, the above
definition is more convenient for our purposes.

The bound \eqref{E:BDBound} 
was improved in several steps by Huxley \cite{Huxley} to $0.4394\ldots$.
In 1988, Maier \cite{Maier} used his matrix method to 
improve the bound to $0.2484\ldots$. 
Recently, the first, third and fourth authors proved 
a best possible result in this direction.

\begin{theorem} \label{T:GPY1}
{\rm (Goldston, Pintz, and Yildirim \cite{GPY})}
\begin{equation*}
 \liminf_{n\to \infty}
    \left( 
       \frac{ p_{n+1}-p_n }{\log p_n} 
      \right) =0.
\end{equation*}
\end{theorem}

The proof of Theorem \ref{T:GPY1} uses, among other things, the 
Bombieri-Vinogradov Theorem.
There are good reasons to believe that the bound in 
\eqref{E:TheoBConclusion} holds for larger values of $Q$.
More formally 
we have the following conjecture.

\begin{HypoB}
Suppose $1/2 < \theta \le 1$.
If $A>0, \epsilon >0$, then
\begin{equation} \label{E:HypoBConclusion}
\sum_{q\le N^{\theta-\epsilon}}  |E^*(N;q,a)| 
\ll_{A,\epsilon} N(\log N)^{-A}.
\end{equation}
\end{HypoB}

If  Hypothesis $BV(\theta)$ is true, 
then we say that the sequence $\varpi$ has 
{\it level of distribution} $\theta$. Thus the Bombieri-Vinogradov Theorem 
shows that $\varpi$ has a level of distribution $1/2.$
The statement that $\varpi  $ has a level of 
distribution $1$ is known as the 
``Elliott-Halberstam Conjecture'' \cite{EH}. Any level of 
distribution larger than $1/2$ will give the following strengthening 
of Theorem \ref{T:GPY1}.

\begin{theorem}\label{T:GPY2}
    {\rm (Goldston, Pintz, and Yildirim \cite{GPY})} 
    If Hypothesis $BV(\theta)$ is true 
    for some $\theta >1/2$, then 
\begin{equation*}
 \liminf_{n\to \infty} (p_{n+1}-p_n) < \infty.
\end{equation*}
If Hypothesis $BV(\theta)$ is true for some $\theta$ with
$4(8-\sqrt{19})/15 = 0.97096\ldots < \theta \le 1$, then 
\begin{equation*}
    \liminf_{n\to \infty} (p_{n+1} -p_n ) \le 16.
\end{equation*}
\end{theorem}

Our first objective here is to give alternative proofs of Theorems
\ref{T:GPY1} and \ref{T:GPY2}. The primary difference in the proofs
here and the proofs in \cite{GPY} comes from the use of Selberg
diagonalization and a different choice of sieve coefficients; this 
will be discussed in more detail below. Our choice of coefficients
allows us to give an elementary treatment of the main terms; we 
will discuss this further after the statement of Theorem \ref{T:Thm6}
below.

Our second objective is to show that the results of \cite{GPY}
can be strengthened if one replaces primes by numbers with 
a fixed number of prime factors. Let $E_k$ denote a number with
 numbers with {\it exactly} $k$ distinct prime factors. This contrasts with
the usual definition of ``almost-prime'', where $P_k$ is used
to denote a number with at most $k$ distinct prime factors. 
Chen \cite{Chen} proved that there are infinitely many primes
$p$ such that $p+2$ is a $P_2.$ While one expects that there
are infinitely many primes $p$ such that $p+2$ is an $E_2$, this 
appears to be as difficult as the twin prime conjecture.
However, we can prove that the limit infimum of
gaps between $E_2$'s is bounded.

\begin{theorem} \label{T:E2Gaps}
Let $q_{n}$ denote the $n^{\rm th}$ number that is a product of exactly 
two primes. Then
\begin{equation*}
\liminf_{n\to\infty} \left( q_{n+1}-q_{n} \right)\le 26.
\end{equation*}
\end{theorem}

The above theorem uses an analogue of the Bombieri-Vinogradov theorem 
for the function $\varpi*\varpi$, which is defined as 
\begin{equation*}
\varpi*\varpi(n) =\sum_{d|n} \varpi(d)\varpi(n/d).
\end{equation*}
Note that $\varpi*\varpi(n)=0$ unless $n$ is a product of two 
primes or $n$ is a square of a prime. 
 
When $(a,r)=1$, we have 
\begin{equation*}
  \sum_{\substack{N < n \le 2N \\ n\equiv a \pmod r}} 
   \varpi*\varpi(n)
   =
   \frac{1}{\phi(r)} \sum_{\chi \pmod r} \bar{\chi}(a) 
       \sum_{N < n \le 2N} \varpi*\varpi(n) \chi(n),
\end{equation*}
and the expected value of this is
\begin{equation}  
 \frac{1}{\phi(r)} 
   \sum_{N < n \le 2N} \varpi*\varpi(n) \chi_0(n),
\end{equation}
where $\chi_0$ is the principal character mod $r$. 
A computation (see Lemma \ref{L:Varpi})
shows that this quantity is asymptotically equal to
\begin{equation} \label{E:R1}
   \frac{N}{\phi(r)}
    \left( \log N + C_0 - 2 \sum_{p|r} \frac{\log p}{p} \right),
\end{equation} 
where $C_0$ is the absolute constant defined in \eqref{E:C0Def}.

Let $E_2(N;r,a)$ be defined by 
\begin{equation*}
\sum_{\substack{N< n\le 2N \\ n \equiv a \pmod r}} \varpi*\varpi(n)= 
     \frac{N}{\phi(r)} 
      \left(\log N + C_0 -2 \sum_{p|r} \frac{\log p}{p} \right) + E_2(N;q,a).
\end{equation*}
In parallel to the definitions of $E(N,q)$ and $E^*(N,q)$, we define
\begin{equation*}
E_2(N,r) =\max_{a, (a,r)=1} |E_2(N;r,a)|, \quad
E^*_2(N,r)= \max_{x\le N} E_2(x,r).
\end{equation*}

\begin{TheoB2} 
For every $A>0$, there exists $B>0$ such that
if $Q\le N^{1/2}\log^{-B}N$
\begin{equation*}
\sum_{r\le Q}  |E_2^*(N,r)| \ll_A N(\log N)^{-A}.
\end{equation*}
\end{TheoB2}

This is a special case of a result of Motohashi \cite{Motohashi}.
Alternatively, one can easily modify the Vaughan's Identity for 
the von Mangoldt function
$\Lambda$ to an identity for $\Lambda*\Lambda$, and then use Vaughan's 
approach (see \cite{Va} or Chapter 28 of \cite{Da}) to the Bombieri-Vinogradov 
Theorem to prove the analogue for $\Lambda*\Lambda$. 
It is then easy to modify this to a result for $\varpi*\varpi$.

We also propose a natural analogue of Hypothesis $BV(\theta)$.

\begin{HypoB2} 
Suppose $1/2 <  \theta \le 1$.
If $A>0, \epsilon >0$, then
\begin{equation} \label{E:HypoB2Conclusion}
\sum_{q\le N^{\theta-\epsilon}}  |E_2^*(N;q)| 
\ll_{A,\epsilon} N(\log N)^{-A}.
\end{equation}
\end{HypoB2} 

From this, we obtain the following conditional result. 

\begin{theorem} \label{T:E2GapsEH}
If Hypotheses $BV(\theta)$ and $BV_2(\theta)$ are 
both true for some $\theta$ with 
$(75-\sqrt{473})/56 = 0.950918 \ldots < \theta \le 1$, then  
\begin{equation*}
    \liminf_{n\to\infty} \left( q_{n+1} - q_{n} \right)\le 6.
\end{equation*}
\end{theorem}

The basic construction for the proofs of Theorems \ref{T:GPY1}
and \ref{T:GPY2} 
was inspired by work of Heath-Brown 
\cite{Heath-Brown} on almost prime-tuples of linear forms. 
Heath-Brown's work was itself a generalization of Selberg's
proof \cite{S1992} that the polynomial $n(n+2)$ will infinitely often have 
at most five prime factors, and in such a way that one of
$n$ and $n+2$ has at most two prime factors, while the other has 
at most three prime factors.

Define
\begin{equation} \label{E:PDefinition}
P(n;\cH)  = \prod_{h\in \cH} (n+h), 
\end{equation}
The central idea is to relate the problem to sums of the form
\begin{equation} \label{E:BasicSum}
\sum_{N< n \le 2N}
 \left( \sum_{d|P(n;\cH)} \lambda_d \right)^2
\end{equation} 
and of the form
\begin{equation} \label{E:BasicSumW}
\sum_{N< n \le 2N} \varpi(n) 
 \left( \sum_{d|P(n;\cH)} \lambda_d \right)^2, 
\end{equation}
where one assumes that $\lambda_d=0$ for $d>R$, and $R$ is a parameter 
that is chosen to control the size of the error term.  
One also assumes that $\lambda_d=0$ when $d$ is not squarefree.

To illustrate the relevance of the sums \eqref{E:BasicSum}
and \eqref{E:BasicSumW}, we discuss one simple application
that is related to the second part of Theorem \ref{T:GPY2}.
Let $\cH$ be an admissible $k$-tuple, and consider the sum
\begin{equation} \label{E:cSDef}
\cS:=\sum_{N< n \le 2N}  
 \left\{  \sum_{h\in \cH} \varpi(n+h)  \,-\, (\log 3N) \right\} 
 \left( \sum_{d|P(n;\cH)} \lambda_d \right)^2.
\end{equation}
For a given $n$, the inner sum is negative unless there are at least 
two values $h_i,h_j\in \cH$ such that $n+h_i,n+h_j$ are primes. 
From Theorems \ref{T:Thm5} and \ref{T:Thm6} below, one can deduce that
if $BV(\theta)$ is true, if $R=N^{\theta-\epsilon}$ for   
$\epsilon>0$, and 
if $0\le \ell \le k$,  then 
\begin{equation*}
    \cS \gtrsim N\fS(\cH)(\log R)^{k+2\ell} (\log N) 
    m(k,\ell,\theta),
\end{equation*}
where
\begin{equation*}
    m(k,\ell,\theta)= 
    \binomial{2\ell}{\ell} \frac{1}{(k+2\ell)!}
    \left\{ 
       \frac{k(2\ell+1)(\theta-\epsilon)}{(k+2\ell+1)(\ell+1)} -1
     \right\}.
\end{equation*}
This last expression is positive, if for example, $k=7$,
$\ell=1$, $\epsilon$ is sufficiently small, and $20/21 < \theta \le 1$. 
Consequently, if $BV(1)$ is true, then for any 
admissible $7$-tuple $\cH$, there are infinitely many $n$ and 
some $h_i,h_j\in \cH$ such that $n+h_i, n+h_j$ are both prime. 
Now 
\begin{equation*}
    \cH=\{11,13, 17,19,23, 29, 31\}
\end{equation*}
is an admissible $7$-tuple. $\cH$ is admissible because if $p\le 7$, 
then
none of the elements in $\cH$ are divisible by $p$, and if
$p>7$, then there are not enough elements in $\cH$ to cover all of the 
residue classes mod $p$. Now any two elements of $\cH$ differ
by at most $20$, so we conclude that if $BV(1)$ is true, then
\begin{equation*}
    \liminf_{n\to \infty} (p_{n+1}-p_n) \le 20.
\end{equation*}
To get the stronger bound of $16$ given in Theorem \ref{T:GPY2} needs
an extra idea; this will be discussed in Section \ref{S:OtherThms}.

The success of the method depends upon making an appropriate choice
for the $\lambda_d$, and this takes us into the realm of the 
Selberg upper bound sieve. 
It is a familiar fact from the theory of this sieve 
that 
\begin{equation*}
   \sum_{\substack{N < n \le 2N \\ d|P(n;\cH)}} 1
   = \frac{N}{f(d)} +r_d,
\end{equation*}
where $f$ is a multiplicative function and $r_d$ is a remainder term. 
(See the first part of Section \ref{S:Thm5} for the formal definition
of $f$.)
Accordingly, an appropriate transformation of the sum in
\eqref{E:BasicSum} leads to consideration of the  bilinear form
\begin{equation} \label{E:BilinearForm1}
    \sum_{d,e} \frac{\lambda_d \lambda_e}{f([d,e])}.
\end{equation}
The typical approach in the Selberg sieve 
is to choose the $\lambda_d$ to minimize the form
in \eqref{E:BilinearForm1}. To make this problem feasible, one needs 
to diagonalize this bilinear form. This can be done 
by making a change of variables
\begin{equation} \label{E:yrDefinition1}
    y_r= \mu(r) f_1(r) \sum_{d} \frac{\lambda_{dr}}{f(dr)},
\end{equation}
where $f_1$ is the multiplicative function defined by $f_1=f*\mu$.
(Note that the sum in \eqref{E:yrDefinition1} is finite because
$\lambda_d=0$ for $d>R$.) The sum in \eqref{E:BilinearForm1} is 
then transformed into
\begin{equation*}
    \sum_{r} \frac{y_r^2}{f_1(r)},
\end{equation*}
and the bilinear form is minimized by taking
\begin{equation} \label{E:Sy}
    y_r = \mu^2(r) \frac{\lambda_1}{V},
\end{equation}
where
\begin{equation*}
    V=\sum_{r < R} \frac{\mu^2(r)}{f_1(r)}.
\end{equation*}
The minimum of the form in \eqref{E:BilinearForm1} is then seen
to be
\begin{equation*}
\frac{\lambda_1^2}{V}.
\end{equation*}
One usually assumes that $\lambda_1=1$, but this is not an essential
element of the Selberg sieve, and it is sometimes useful to assign 
some other nonzero value to $\lambda_1$.

The sum in \eqref{E:BasicSumW} can be treated in a similar way. 
However, the corresponding function $f$ must be replaced by a slightly
different function $f^*$, which will be defined in Section \ref{S:Thm6}.
Therefore, the optimal choice of $\lambda_d$ is different from the 
optimal choice for the sum in \eqref{E:BasicSum}. However, the basic
structure of our approach requires that the same choice of $\lambda_d$
be used for both sums. We therefore face the problem of making a 
choice of $\lambda_d$ that works reasonably well for both problems. 
A similar choice was faced by Selberg and Heath-Brown in their 
earlier mentioned work, and they made this choice in different ways.
Selberg \cite{S1992} made a choice of $\lambda_d$ that was optimal for 
one problem, 
and was able to successfully analyze the effect of this choice for the 
other problem. Heath-Brown \cite{Heath-Brown} chose
\begin{equation*}
    \lambda_d = 
     \begin{cases}
	\displaystyle \mu(d) \left(\frac{\log R/d}{\log R}\right)^{k+1} &
	  \text{ if $d< R$,}\\
	  0 &\text{ otherwise;}
      \end{cases}
\end{equation*}
$k$ being the number of linear forms under consideration. 
While this choice is not optimal for either problem, it is asymptotically 
optimal for both problems.

Inspired by Heath-Brown's choice, Goldston, Pintz, and Yildirim 
\cite{GPY} chose
\begin{equation}  \label{E:GPYChoice}
    \lambda_{d,\ell}= 
      \begin{cases}
	  \displaystyle \mu(d) \frac{(\log R/d)^{k+\ell}}{(k+\ell)!} &
	   \text{ if $d< R$,}\\
	  0 &\text{ otherwise.}
      \end{cases}
\end{equation}
Here, $\ell$ is a non-negative integer to be chosen in due course,
with $\ell\le k$.
With the exponent $k+\ell$, one is effectively using  
a $k+\ell$-dimensional sieve on a $k$-dimensional
sieve problem. In an upper bound sieve, it is optimal to take the 
dimension of the sieve to be the same as the dimension of the problem.
In the problems considered here, however, it is not the upper bound 
but the ratio of the quantities  in \eqref{E:BasicSum} and
\eqref{E:BasicSumW} that is relevant. The presence of the parameter 
$\ell$ is essential for the success of their method.  

In the current exposition, we make a choice that is a hybrid of 
the above and Selberg's original approach. Our choice is most 
easily described in terms of $y_r$. We choose
\begin{equation} \label{E:yrChoice}
    y_{r,\ell}=y_{r,\ell}(\cH)=
     \begin{cases}
	 \displaystyle \frac{ \mu^2(r) \fS(\cH) (\log R/r){^\ell}}{\ell!} &
	  \text{ if $r< R$,}\\
	  0 &\text{ otherwise.}
      \end{cases}
\end{equation}
As motivation for this choice,
we note that $y_{r,0}$ is the optimal choice given
in \eqref{E:Sy} with $\lambda_1=V\fS(\cH)$. 
Moreover, one can show that 
\begin{equation*}
    \mu(r)f_1(r) \sum_{d<R/r} 
       \frac{\mu(dr)}{f(dr)} 
         \frac{\log^{k+\ell} (R/rd)}{(k+\ell)!}
   \asymptotic \frac{ \fS(\cH) (\log R/r)^\ell}{\ell!}
\end{equation*}
when $r$ is not too close to $R$.
In other words, the choice of $\lambda_{d,\ell}$ in \eqref{E:GPYChoice} gives
a value of $y_r$ that is asymptotic to the expression in 
\eqref{E:yrChoice}.

One can use 
\eqref{E:yrDefinition1} and  M\"obius inversion
to deduce that
\begin{equation} \label{E:lambda-yRelation}
    \frac{\lambda_{d,\ell}}{f(d)} =
     \mu(d) \sum_r \frac{y_{dr,\ell}}{f_1(rd)},
\end{equation}
and so, when the choice of $y_{r,\ell}$ of \eqref{E:yrChoice} is specified,
one obtains
\begin{equation} \label{E:lambdaChoice}
\lambda_{d,\ell}  =
\mu(d) \frac{f(d)}{f_1(d)} \frac{\fS(\cH)}{\ell!} 
\sum_{\substack{r< R/d\\(r,d)=1}}
\frac{ \mu^2(r) }{f_1(r)}  (\log R/rd)^{\ell} 
\end{equation}
when $d < R$.
With this choice of $\lambda_{d,\ell}$, we set
\begin{equation} \label{E:LambdaRChoice}
\Lambda_R(n;\cH,\ell)=\sum_{d|P(n;\cH)} \lambda_{d,\ell}.
\end{equation}
As we shall see, this choice $\lambda_{d,\ell}$ allows us to give
elementary estimates for the main terms in 
\eqref{E:BasicSum} and \eqref{E:BasicSumW}.

We also define 
\begin{equation} \label{E:betaDef}
    \beta(\cH)= \sum_{p} \frac{(k-\nu_p(\cH))\log p}{p}.
\end{equation}
This sum is finite because $\nu_p=k$ for sufficiently large $p$.

Theorems \ref{T:GPY1} through \ref{T:E2GapsEH} 
will be derived fairly easily from the 
following results.

\begin{theorem}\label{T:Thm5} Suppose that 
$\cH=\{h_1,\ldots,h_k\}$ is an admissible set, and that 
$0\le \ell_1, \ell_2 \le k$.
If $R\le N^{1/2-\epsilon}$ then
\begin{align} \label{E:Thm5}
\sum_{N < n \le 2N} & 
\Lambda_{R} (n;\cH,\ell_1) \Lambda_{R}(n;\cH,\ell_2)= \\
& \binom{\ell_1+\ell_2}{\ell_1} \fS(\cH) 
N\frac{(\log R)^{k+\ell_1+\ell_2}}{(k+\ell_1+\ell_2)!}
\SE. \notag 
\end{align}
The implied constant depends at most on $k$. 
\end{theorem}

\begin{theorem} \label{T:Thm6}
Suppose that $\cH=\{h_1,\ldots, h_k\}.$
Suppose further that Hypothesis $BV(\theta)$
is true and $R\le N^{(\theta-\epsilon)/2}$. 
If $h_0\in \cH$, $\cH$ is admissible, and
$0\le \ell_1,\ell_2\le k$, then 
\begin{align} \label{E:Thm6Part1}
\sum_{N < n \le 2N} & \varpi(n+h_0) 
\Lambda_{R}(n;\cH,\ell_1) \Lambda_{R}(n;\cH,\ell_2)  = \\
   & \binom{\ell_1+\ell_2+2}{\ell_1+1} N \fS(\cH) 
         \displaystyle{ 
               \frac{(\log R)^{k+\ell_1+\ell_2+1} }{(k+\ell_1+\ell_2+1)!}  
                   } \SE\notag.
\end{align}

If $h_0\notin \cH$, $\cH^{0}=\cH\union\{h_0\}$ is admissible, 
and $1\le \ell_1,\ell_2 \le k$ then 
\begin{align} \label{E:Thm6Part2}
\sum_{N < n \le 2N} & \varpi(n+h_0) 
\Lambda_{R}(n;\cH,\ell_1) \Lambda_{R}(n;\cH,\ell_2)  = \\
   & \binom{\ell_1+\ell_2}{\ell_1} N \fS(\cH^0) 
         \displaystyle{ 
               \frac{(\log R)^{k+\ell_1+\ell_2} }{(k+\ell_1+\ell_2)!}  
                   }
\left\{ 1+O(\beta(\cH^0)\fS(\cH^0)/\log R \right\}. \notag
\end{align}
The implied constants depend at most on $k$.
\end{theorem}
 
With a bit more work, we could allow $\ell_1$ or $\ell_2$ to be $0$
in \eqref{E:Thm6Part2}. However, we omit this because
the only place we use this 
result is in the proof of Theorem \ref{T:GPY1}, 
where we will have $\ell_1=\ell_2 >0$.

Analogues of Theorems \ref{T:Thm5} and \ref{T:Thm6} are given in 
\cite{GPY} for $\lambda_{d,\ell}$  given by \eqref{E:GPYChoice}. 
The corresponding main terms in \cite{GPY} are evaluated with the help 
of contour integrals in two variables and zero-free regions for the 
Riemann-zeta function. On the other hand, with the choice of 
$\lambda_{d,\ell}$ given
in \eqref{E:lambdaChoice}, we are able to give an elementary treatment
of the main terms in Theorems \ref{T:Thm5} and \ref{T:Thm6}.

\begin{theorem} \label{T:Thm7}
Suppose that 
$\cH=\{h_1,\ldots,h_k\}$ is an admissible set, and that 
$0\le \ell_1, \ell_2 \le k$.
Suppose that Hypotheses $BV(\theta)$ and $BV_2(\theta)$ are 
both satisfied, and 
$R\le N^{(\theta-\epsilon)/2}$. If $h_0\in \cH$,  then
\begin{align*}
\sum_{N < n \le 2N} & \varpi*\varpi(n+h_0) 
      \Lambda_{R}(n;\cH,\ell_1)  \Lambda_{R}(n;\cH,\ell_2) = \\
    &  \left\{
       \binom{\ell_1+\ell_2+2}{\ell_1+1}  (N\log N)\fS(\cH) 
       \frac{ (\log R)^{k+\ell_1+\ell_2+1}} {(k+\ell_1+\ell_2+1)!} 
          \right.\\   
    & \left.+2T(k,\ell_1,\ell_2) N \fS(\cH) 
        \frac{(\log R)^{k+\ell_1+\ell_2+2}} 
               {(k+\ell_1+\ell_2+2)!}  \right\} 
           \SE,
\end{align*}
where
\begin{equation*}
    T(k,\ell_1,\ell_2) = 
           -\binom{\ell_1+\ell_2+3}{\ell_2+1}
            - \binom{\ell_1+\ell_2+3}{\ell_1+1}
	     + \binom{\ell_1+\ell_2+2}{\ell_1+1}.
\end{equation*}
The implied constant depends at most on $k$.
\end{theorem}

The reader will note that the sums considered here are more 
general than the sums in \eqref{E:BasicSum} and 
\eqref{E:BasicSumW}--the latter correspond to the case
$\ell_1=\ell_2=\ell$.
We will see in Section \ref{S:OtherThms} that
this extra flexibility is useful in applications.

We also remark that the proof of Theorem \ref{T:GPY1} requires 
averaging over a set of $\cH$, where the elements of $\cH$ can
be as large as $\log R$. Accordingly, we shall take some extra effort 
to make our estimates uniform in $h$ under the assumption that $h\le 
\log N$. For our results, it is not necessary to make the estimates in
Theorems \ref{T:Thm5} through \ref{T:Thm7} uniform in $k$.

The implied constants in the error terms of Theorems 6 and 7 are
ineffective due to the use of the Bombieri-Vinogradov Theorem, which
uses the Siegel-Walfisz Theorem. However, the constants can be made 
effective by using the procedure of Section 13 of \cite{GPY}. 
This procedure deletes the greatest prime factor of the eventually 
existing exceptional modulus from the sieve process.

The paper \cite{GPY} gives an unconditional proof of a quantitative
version of Theorem \ref{T:GPY1}; that
\begin{equation} \label{T:GPY3}
\liminf_{n\to\infty}
\frac {(p_{n+1}-p_n)}
       {\log p_n (\log\log p_n)^{-1} \log\log\log\log p_n} 
< \infty,
\end{equation}
and this result requires that the  estimates in Theorems
\ref{T:Thm5} and \ref{T:Thm6} be uniform in $k$. 
In a forthcoming paper, Goldston, Pintz, and Yildirim will improve
\eqref{T:GPY3} to 
\begin{equation}  
\liminf_{n\to\infty}
\frac {(p_{n+1}-p_n)}
       {(\log p_n)^{1/2} (\log\log p_n)^2} 
< \infty.
\end{equation}

The function $\varpi*\varpi$ used in Theorem \ref{T:Thm7}
is convenient for calculations, but it is not optimal for 
applications. In a future paper we will show that by using other 
functions supported on $E_2$'s, the bound
in Theorem \ref{T:E2Gaps} can be improved to 8 and the 
allowable range for $\theta$ in Theorem \ref{T:E2GapsEH}
can be improved to $0.51 < \theta \le 1$. We will also
show that  there is a constant $C$ such that for any 
positive integer $r$,
$$\liminf_{n\to\infty} (q_{n+r}-q_n) \le C r  e^r.$$

{\it Notation:}
The letters $R,N$ denote real variables tending to infinity.
The letter $p$ is always used to denote a prime. The letters 
$d,e,r$ are usually squarefree numbers; the letters $m,n$ are usually 
positive integers.
The notation $\omega(n)$ is used to denote the number of distinct prime 
factors of $n$. We use $\rho$ to denote the function
\begin{equation*}
    \rho(r)= 1+ \sum_{p|r} \frac{\log p}{p}.
\end{equation*}
The letters $S,\cL,U,$ and $V$,  with or without subscripts, are often 
used to denote sums. The meanings of these symbols are local to sections;
e.g., the meaning of $S_1$ in Section \ref{T:Thm6} is different from 
the meaning of $S_1$ in Section \ref{T:Thm7}.

We use  $\sumflat_\null$ to denote a summation over squarefree integers.
In general, the constants implied by ``$O$'' and ``$\ll$'' will 
depend on $k$. Any other dependencies will be explicitly noted.
As noted before, $k$ is the size of $\cH$; we always assume that 
$k\ge 2$. The parameter $\ell$, with or without subscript,
is an integer with $0\le \ell\le k$.

%This work is preliminary and subject to further checking.
%This manuscript is intended for limited distribution.

\section{Preliminary Lemmas} \label{S:Lemmas}

The following two lemmas are classical estimates that have proved 
useful for handling remainder terms that arise in the Selberg sieve. 
The results can be found in Halberstam and Richert's 
book (\cite{HR}, Lemmas 3.4 and 3.5). We reproduce the proofs here
since they are quite short.

\begin{lemma}  \label{L:HR34}
    For any natural number $h$ and for $x\ge 1$, 
    \begin{align*}
	\sumflat_{d\le x} \frac{h^{\omega(d)}}{d} 
	   \le & (\log x +1)^h, \\
	 \sumflat_{d\le x} {h^{\omega(d)}}
	  \le & x(\log x+ 1)^h.
    \end{align*}
\end{lemma}

\begin{proof} For the first inequality, we note that the sum on the 
left is
\begin{equation*}
     \sum_{d_1\ldots d_h \le x} 
	\frac{\mu^2(d_1\ldots d_h)}{d_1\ldots d_h}
	\le \left( \sum_{n\le x} \frac{1}{n} \right)^h \le (\log x+1)^h.
\end{equation*}
For the second inequality, we note that the left-hand side is
at most 
\begin{equation*}
     x \sumflat_{d\le x} \frac{h^{\omega(d)}}{d},
\end{equation*}
and we appeal to the first inequality.
\end{proof}

\begin{lemma} \label{L:HR35}
Assume Hypothesis $BV(\theta)$, and let $h$ be a positive integer.
Given any positive constant $U$ and any $\epsilon >0$,
then
\begin{equation*}
    \sumflat_{d< N^{\theta-\epsilon}} h^{\omega(d)} E^*(N,d)
       \ll_{U,h,\epsilon} N(\log  N)^{-U}.
\end{equation*}

Similarly, if Hypothesis $BV_2(\theta)$ is assumed,  then 
\begin{equation*}
    \sumflat_{d< N^{\theta-\epsilon}} h^{\omega(d)} E_2^*(N,d)
       \ll_{U,h,\epsilon}  N(\log  N)^{-U}.     
\end{equation*}
\end{lemma}

\begin{proof}
    We begin by noting the trivial estimate $E^*(N,d) \ll N(\log N)/d$.
By Cauchy's inequality
\begin{align*}
    \sumflat_{d<N^{\theta-\epsilon}} h^{\omega(d)} E^*(N,d) 
    \le & \left( N\log N \sumflat_{d<N^{\theta-\epsilon}} 
		   \frac{h^{2\omega(d)}}{d} \right)^{1/2}
	\left( \sumflat_{d< N^{\theta-\epsilon}} E^*(N,d) \right)^{1/2} 
	\\
	& \ll_{h,\epsilon,A} N (\log N)^{(h^2-A+1)/2}.
\end{align*}
We have used Lemma \ref{L:HR34} and Hypothesis $BV(\theta)$ in the last line.
The first result follows by taking $A=h^2+1+2U$. The second result is
proved similarly; one uses the trivial bound 
$E_2^*(N,d) \ll N(\log N)^2/d$.
\end{proof}

 \begin{lemma} \label{L:BetaIntegral}
    If $a,b$ are positive real numbers, both at least $1$,
then 
\begin{equation*}
\int_1^x (\log x/u)^{a-1} (\log u)^{b-1} \frac{du}{u}
= (\log x)^{a+b-1} \frac{\Gamma(a)\Gamma(b)}{\Gamma(a+b)}.
\end{equation*}
\end{lemma}
\begin{proof} Upon making the change of variables $u=x^v$,
    the left-hand side becomes
    \begin{equation*}
	(\log x)^{a+b-1} \int_0^1 (1-v)^{a-1} v^{b-1} dv.
    \end{equation*}
The result follows by the standard formula for the beta-integral.
\end{proof}

Our next lemma is another standard result in the theory of sieves.

\begin{lemma} \label{L:Wirsing}
Suppose that $\gamma$ is a multiplicative function,
and suppose that there positive real numbers 
$\kappa,A_1,A_2, L$ such that
\begin{equation} \label{E:Omega1}
	0\le \frac{\gamma(p)}{p} \le 1-\frac{1}{A_1},
\end{equation}
and 
\begin{equation} \label{E:Omega2}
	-L \le \sum_{w\le p < z} 
	         \frac{\gamma(p) \log p}{p}
		 -\kappa \log \frac{z}{w} 
		\le A_2
\end{equation}
if $2 \le w \le z$. 
Let $g$ be the multiplicative function defined by 
\begin{equation} \label{E:gDefinition}
    g(d) = \prod_{p|d} \frac{\gamma(p)}{p-\gamma(p)}.
\end{equation}
Then 
\begin{equation*}
    \sumflat_{d < z} g(d) =
     c_{\gamma} \frac{(\log z)^{\kappa}}{\Gamma(\kappa+1)}
        \left\{ 1+O_{A_1,A_2,\kappa}\left( \frac{L}{\log z} \right) \right\},
\end{equation*} 
where
\begin{equation*}
    c_{\gamma} = 
     \prod_p \left( 1-\frac{\gamma(p)}{p} \right)^{-1}
             \left( 1-\frac{1}{p} \right)^{\kappa}.
\end{equation*}
\end{lemma}

This is a combination of Lemmas 5.3 and 5.4 of Halberstam and 
Richert's book \cite{HR}. 
In \cite{HR}, the hypothesis \eqref{E:Omega1} is denoted
$(\Omega_1)$, and hypothesis \eqref{E:Omega2} is denoted
$(\Omega_2(\kappa,L))$.
As indicated above, the constant implied by ``$O$'' may depend on 
$A_1,A_2,\kappa$, but it is independent of $L$. 
This will be important in our applications.

\begin{lemma} \label{L:ExtendedWirsing}
Suppose that $\gamma$ and $g$ satisfy the same hypotheses
as in the previous lemma. If $a$ is a non-negative integer, then
\begin{equation*}
    \sumflat_{r< R} g(r) (\log R/r)^a = 
     c_{\gamma} \frac{\Gamma(a+1)}{\Gamma(\kappa+a+1)} (\log R)^{\kappa+a}
        +O_{A_1,A_2,\kappa,a}\left(L (\log R)^{\kappa+a-1}\right).
\end{equation*}
\end{lemma}

\begin{proof} When $a=0$, this is Lemma \ref{L:Wirsing}.
If $a > 0$, then 
\begin{align*}
    \sumflat_{r<R} g(r) (\log R/r)^a & =
       a\sumflat_{r<R} g(r) \int_r^R (\log R/z)^{a-1} \frac{dz}{z} \\
    & = \int_1^R \frac{a (\log R/z)^{a-1}}{z} \sumflat_{r < z} g(r) dz.
\end{align*}
Using Lemma \ref{L:Wirsing}, we see that the above is
\begin{align*}
    \int_1^R & \frac{a (\log R/z)^{a-1}}{z} 
         \left\{
           \frac{ c_{\gamma}(\log z)^{\kappa}}{\Gamma(\kappa+1)}
	     + O(L(\log z)^{\kappa-1} ) \right\} \, dz\\
     = & \frac{ac_{\gamma}}{\Gamma(\kappa+1)}
            \int_1^R (\log R/z)^{a-1} (\log z)^{\kappa} \frac{dz}{z} 
	     +O\left( 
                 aL \int_1^R (\log R/z )^{a-1} (\log z)^{\kappa-1} \frac{dz}{z}
	        \right).
\end{align*}
The desired result follows from using Lemma \ref{L:BetaIntegral}.
\end{proof}

\begin{lemma} \label{L:HBound}If $\cH$ is admissible and $|h_i|\le h$
   for all $h_i\in \cH$,  then
\begin{equation} \label{E:BetaIneq}
        1 \ll   \beta(\cH) \ll \log\log 10h
\end{equation}
and there is a constant $b_k$ (depending only on $k$) such that
\begin{equation} \label{E:fSIneq}
          \fS(\cH) \ll (\log \log 10h)^{b_k}.
\end{equation}
\end{lemma}

\begin{proof} Without loss of generality, we 
   may assume that $h\ge 100$; this  will simplify the writing of 
   logarithms.
We note that $\nu_p< k$ if and only if $p|\Delta(\cH)$, where
\begin{equation} \label{E:DeltaDef}
	\Delta=\Delta(\cH):=\prod_{1\le i< j \le k} |h_i-h_j|.
\end{equation}
Therefore
\begin{equation*}
    \beta(\cH)=\sum_{p|\Delta} (k-\nu_p) \frac{\log p}{p},
\end{equation*}
where we have written $\nu_p$ as an abbreviation for $\nu_p(\cH)$.
We may assume without loss of generality that 
$\Delta \ge 100$.

Now $\nu_2=1$ whenever $\cH$ is admissible, so we see that
$\beta(\cH) \ge \log 2/2$. In the opposite direction,
we have
\begin{align*}
    \beta(\cH) & \ll 
     \sum_{p\le \log \Delta} \frac{\log p}{p}
     + \sum_{\substack{p|\Delta \\ p > \log \Delta}}
              \frac{\log\log \Delta}{\log \Delta} \\
	      & \ll \log\log \Delta +
	            \frac{\log\log \Delta}{\log \Delta}
		     \frac{\log \Delta}{\log\log \Delta} \\
	      & \ll \log\log \Delta + 1.
\end{align*}
Finally, note that $\Delta\le h^{k^2}$, so that 
$\log \Delta \ll \log h$. This completes the proof of 
\eqref{E:BetaIneq}.

Now consider $\fS(\cH)$. 
 From the definition of $\fS(\cH)$, we see that
\begin{equation*}
	\log \fS(\cH) = \sum_{p} 
	  \left\{
	   \left(\frac{k-\nu_p}{p} \right) + 
	          O\left( \frac{1}{p^2} \right)
           \right\}
	   \ll 1+ \sum_{p|\Delta} \frac{1}{p}.
\end{equation*}
The last sum may be bounded in a manner similar to that used 
for $\beta(\cH)$. We have
\begin{align*}
    \sum_{p|\Delta} \frac{1}{p} & \le 
       \sum_{p\le \log \Delta} \frac{1}{p} +
        \sum_{\substack{p|\Delta \\ p > \log \Delta }}
	     \frac{1}{\log \Delta} \\
	 & \ll \log\log\log \Delta + 
	     \frac{1}{\log \Delta} \frac{\log \Delta}{\log\log \Delta} \\
	 & \ll \log\log\log \Delta.
\end{align*}
As noted before, $\log \Delta \ll \log h$. Therefore, there is some 
constant $b_k$ such that $\log \fS(\cH) \le b_k \log \log \log h$,
and \eqref{E:fSIneq} follows.
\end{proof}

In our final lemma of this section, we give a computation that 
was used in \eqref{E:R1}.

\begin{lemma} \label{L:Varpi}
    Suppose that $q$ is an integer with all of its
    prime divisors less than $\sqrt N$. Then there is some 
    absolute constant $c$ such that
    \begin{equation*}
	\sum_{\substack{N < n \le 2N\\(n,q)=1}} \varpi*\varpi(n)
	= 2N \left(\log N+ C_0 - \sum_{p|q} \frac{\log p}{p} \right)
	 + O(N\exp(-c\sqrt{\log N})),
    \end{equation*}
    where
    \begin{equation} \label{E:C0Def}
	C_0= 2\log 2- 2\gamma -1 - 2\sum_{p} \frac{\log p}{p(p-1}.
     \end{equation}
\end{lemma}
\begin{proof}
We first use the hyperbola method to write
\begin{align*}
    \sum_{n\le x} \varpi*\varpi(n) = &
    2\sum_{m\le \sqrt x} \varpi(m) \sum_{n\le x/m} \varpi(n) - 
    \left( \sum_{m\le \sqrt x} \varpi(m) \right)^2 \\
     = & 2 x \sum_{p\le \sqrt x} \frac{\log p}{p}  - x + 
        O\left(x\exp(-c\sqrt{\log x})\right).
\end{align*}
Next, we use the classical estimate
\begin{equation*}
    \sum_{p\le x} \frac{\log p}{p} =
     \log x -\gamma - \sum_p \frac{\log p}{p(p-1)} + 
     O(\exp(-c\sqrt{\log x})),
\end{equation*}
to get 
\begin{equation} \label{E:Varpi1}
   \sum_{n\le x} \varpi*\varpi(n) = 
    x\log x + C_1 x + O(x\exp(-c\sqrt{\log x})),
\end{equation}
where
\begin{equation*}
    C_1= -2\gamma-2\sum_{p}\frac{\log p}{p(p-1)} - 1.
\end{equation*}
We use \eqref{E:Varpi1} with $x=N$, $x=2N$, and take 
differences to get  
\begin{equation} \label{E:Varpi2}
    \sum_{N < n \le 2N} \varpi*\varpi(n)
    = N\log N + N C_0 + O(N\exp(-c\sqrt{\log N})).
\end{equation}
Finally, we note that for a given integer $q < \sqrt{N}$,
\begin{align}
   \sum_{p|q} \sum_{\substack{N< n \le 2N\\(n,q)=p}} 
   \varpi*\varpi(n) 
   =  & 2\sum_{p|q} \log p \sum_{N/p< n \le 2N/p} \varpi(n) 
    \label{E:Varpi3} \\
   = & 2N\sum_{p|q} \frac{\log p}{p} + O(N\exp(-c\sqrt{\log N})). 
   \notag
\end{align}
The lemma follows by combining \eqref{E:Varpi2} and \eqref{E:Varpi3}.
\end{proof}

\section{Proof of Theorem \ref{T:Thm5}} \label{S:Thm5}

As we noted in the introduction, we take 
$\nu_p(\cH)$ to be 
the number of distinct residue classes mod $p$ in $\cH$.
We extend this definition to arbitrary squarefree moduli $d$ as follows.
Let $\bZ_d$ be the ring of integers 
mod $d$ and define
\begin{equation} \label{E:OmegaDef}
    \Omega_d(\cH) = 
    \{ a\in \bZ_d : P(a;\cH) \equiv  0 \pmod d\},
\end{equation}
We define $\nu_d(\cH)$ to be the cardinality of $\Omega_d(\cH)$.

Assume that $d_1,d_2$ are squarefree numbers with $(d_1,d_2)=1$. 
The Chinese Remainder Theorem gives an isomorphism 
\begin{equation} \label{E:CRT}
    \xi: \bZ_{d_1} \times \bZ_{d_2} \to \bZ_{d_1d_2}.
\end{equation}
The set $\Omega_{d_1d_2}(\cH)$ is the image of 
$\Omega_{d_1}(\cH) \times \Omega_{d_2}(\cH)$ 
under the isomorphism  $\xi$, so $\nu_d(\cH)$ is multiplicative. 

Throughout this section, we will take $\cH$ to be a fixed admissible
set, and we will usually write $\nu_d$ in place of $\nu_d(\cH)$.

The left-hand side of \eqref{E:Thm5} is
\begin{align} \label{E:Thm1Step1}
\sum_{N < n \le 2N} &
    \left( \sum_{d|P(n;\cH)} \lambda_{d,\ell_1}  \right)
    \left( \sum_{e|P(n;\cH)} \lambda_{e,\ell_2}  \right)  \\
& = \sum_{d,e} \lambda_{d,\ell_1} \lambda_{e,\ell_2}
       \sum_{\substack{N < n \le 2N \\ [d,e]|P(n;\cH) }} 1 \notag \\
& = N \sum_{d,e} \frac{\lambda_{d,\ell_1} \lambda_{e,\ell_2}}{f([d,e])} 
  + O\left(
       \sum_{d,e} |\lambda_{d,\ell_1}\lambda_{e,\ell_2} r_{[d,e]}|\right)  \notag \\
& = NS_1 + O(S_2),\notag
\end{align} 
say, where
\begin{equation} \label{E:fDef-Sec2}
f(d) = \frac{d}{\nu_d},
\end{equation}
and
\begin{equation*}
    r_d= \sum_{\substack{ N< n\le 2N \\ d|P(n;\cH)}} 1 - \frac{N}{f(d)}.
\end{equation*}

The estimates of $S_1$ and $S_2$ require the following two lemmas.
\begin{lemma} \label{L:ExtendedWirsing-f}
We have
    \begin{equation*}
	\sum_{r< R} \frac{\mu^2(r)}{f_1(r)} (\log R/r)^\ell
	= \frac{\ell! (\log R)^{k+\ell} }{\fS(\cH)(k+\ell)!}\SE.
\end{equation*}
\end{lemma}

\begin{proof} 
We apply Lemma \ref{L:ExtendedWirsing} with
\begin{equation*}
    \gamma(p)=\nu_p, \quad
    g(p)=\frac{\nu_p}{p-\nu_p} = \frac{1}{f_1(p)}.
\end{equation*}
Now $\nu_p\le \min(k,p-1),$ so \eqref{E:Omega1} holds with $A_1=k+1$.
Moreover, 
\begin{equation*}
    -\beta(\cH)  \le \sum_{w \le p < z} \frac{(\nu_p-k)\log p}{p}  \le 0
\end{equation*}
and 
\begin{equation*}
	\sum_{w \le p < z } \frac{\log p}{p} = \log(z/w) + O(1).
\end{equation*}
Therefore \eqref{E:Omega2} holds with $\kappa=k$,
    $A_2$ some  constant depending only on $k$, and 
\begin{equation*}
    L \ll 1+\beta(\cH) \ll \beta(\cH).
\end{equation*}
Finally, we note that
    \begin{equation*}
    c_{\gamma} = \prod_p 
		    \left( 1-\frac{\nu_p}{p} \right)^{-1}
		    \left( 1 -\frac{1}{p} \right)^k
		= \frac{1}{\fS(\cH)}.
    \end{equation*}        
\end{proof}

\begin{lemma} \label{L:lambdaBound}
Let $\lambda_{d,\ell}$ be as defined in 
\eqref{E:lambdaChoice}. If $d <R$ and $d$ is squarefree, then
\begin{equation*}
    |\lambda_{d,\ell}| \ll (\log R)^{k+\ell}.
\end{equation*}
\end{lemma}

\begin{proof}    
From \eqref{E:lambdaChoice}, we see that if $d$ satisfies the 
hypotheses of the lemma, then
\begin{align*}
    |\lambda_{d,\ell}| &  =
      \frac{\fS(\cH)}{\ell!} \frac{f(d)}{f_1(d)} 
       \sum_{\substack{r< R/d \\ (r,d)=1}} 
	  \frac{\mu^2(r)}{f_1(r)} (\log R/rd)^\ell \\
	& = \frac{\fS(\cH)}{\ell!} \sum_{t|d} \frac{1}{f_1(t)} 
		 \sum_{\substack{r< R/d \\ (r,d)=1}} 
		    \frac{\mu^2(r)}{f_1(r)} (\log R/rd)^\ell. 
\end{align*}
We move the factor $1/f_1(t)$ inside the sum and write $s=rt$ to get
\begin{align*}
    |\lambda_{d,\ell}| &  =
    \frac{\fS(\cH)}{\ell!} \sum_{t|d}   
		     \sum_{\substack{r< R/d \\ (r,d)=1}} 
			\frac{\mu^2(r)}{f_1(rt)} (\log R/rd)^\ell \\
   & = \frac{\fS(\cH)}{\ell!} \sum_{t|d}   
		     \sum_{\substack{s< Rt/d \\ (s,d)=t}} 
			\frac{\mu^2(s)}{f_1(s)} (\log Rt/sd)^\ell.
\end{align*}
For any $t|d$, we have $Rt/d<R$, so 
\begin{equation*}
   |\lambda_{d,\ell}|  \le \frac{\fS(\cH)}{\ell!}  (\log R)^{\ell}
       \sum_{t|d}
	\sum_{\substack{s < R\\(s,d)=t}} \frac{\mu^2(s)}{f_1(s)}.
\end{equation*}
Now for any $s<R$, there is a unique $t|d$ such that $(s,d)=t$.
Therefore
\begin{equation*}
    |\lambda_{d,\ell}|  \le \frac{\fS(\cH)}{\ell!}  (\log R)^{\ell}
            \sum_{s<R} \frac{\mu^2(s)}{f_1(s)}.
\end{equation*}

To complete the proof, we observe that
\begin{align*}
    \sum_{s<R} \frac{\mu^2(s)}{f_1(s)} & 
       \le \prod_{p < R}  \left( 1+ \frac{1}{f_1(p)}\right) \\
      & =  \prod_{p < R} \left( 1-\frac{\nu_p}{p} \right)^{-1}
                          \left( 1-\frac{1}{p} \right)^{k}
	   \prod_{p < R} \left(1-\frac{1}{p} \right)^{-k} \\
      & \ll \frac{(\log R)^k}{\fS(\cH)}.
\end{align*}
\end{proof}

We now treat $S_1$ and $S_2$. For $S_1$, we begin by writing
\begin{align*}
    S_1 = & \sum_{d,e} 
	     \frac{\lambda_{d,\ell_1} \lambda_{e,\ell_2}}{f(d)f(e)}
		\sum_{\substack{r|d\\r|e}} f_1(r) \\
	=  & \sumflat_r f_1(r) 
		\left( \sum_{d} \frac{\lambda_{dr,\ell_1}}{f(dr)} \right)
		\left( \sum_{e} \frac{\lambda_{er,\ell_1}}{f(er)} \right) \\
	= & \sumflat_r \frac{y_{r,\ell_1} y_{r,\ell_2}}{f_1(r)} \\
	= & \frac{\fS(\cH)^2}{\ell_1! \ell_2!} 
		 \sum_{r < R} 
		  \frac{\mu^2(r) \log^{\ell_1+\ell_2} (R/r)}{f_1(r)}.
\end{align*}

Lemma \ref{L:ExtendedWirsing-f} now yields the estimate
\begin{equation*}
    S_1 = \binomial{\ell_1+\ell_2}{\ell_1} \fS(\cH) 
	   \frac{(\log R)^{k+\ell_1+\ell_2}}{(k+\ell_1+\ell_2)!}\SE.
\end{equation*}

For $S_2$, we first note that
\begin{equation*}
    |r_d| \le \nu_d \le k^{\omega(d)}.
\end{equation*} 
We also have the bound for $\lambda_{d,\ell}$ given in Lemma 
\ref{L:lambdaBound}.
Therefore  
\begin{align*} 
    S_2  & =
       \sum_{d,e< R}|\lambda_{d,\ell_1}\lambda_{e,\ell_2}r_{[d,e]}| 
    \\ 
	 & \ll (\log R)^{2k+\ell_1+\ell_2}  \sumflat_{d,e\le R} k^{\omega([d,e])} 
    \notag\\
	 &  \ll (\log  R)^{4k}  \sumflat_{r < R^2} (3k)^{\omega(r)}.
    \notag 
\end{align*}
Using Lemma 1, we get
\begin{equation} \label{E:Thm1S2}
     S_2 \ll R^2 (\log R)^{7k}  
       \ll (N/\log N)  
\end{equation}
provided $R< N^{1/2-\epsilon}$.

Theorem \ref{T:Thm5} follows by combining the above estimates for 
$S_1$ and $S_2$.

\section{Proof of Theorem \ref{T:Thm6}, Part 1}\label{S:Thm6}
In this section, we consider Theorem 6 under the
assumption that  $h_0\in \cH$.  
Our problem is translation invariant in $\cH$, so we may,
without loss of generality, assume that $h_0=0$ and $0\in \cH$.

Let $\cL$ denote the sum on the left-hand side of \eqref{E:Thm6Part1}.
Then
\begin{equation} \label{E:Thm2LDef}
\cL=\sum_{d,e} \lambda_{d,\ell_1} \lambda_{e,\ell_2} 
 \sum_{\substack{N < n \le 2N\\ [d,e]|P(n;\cH)}} \varpi(n)  =
\sum_{d,e} \lambda_{d,\ell_1} \lambda_{e,\ell_2}
  \sum_{a\in \Omega_{[d,e]}(\cH)} 
  \sum_{\substack{N < p \le 2N \\ p\equiv a \pmod{[d,e]} }}
  \log p.
\end{equation}

Now all prime divisors of $[d,e]$ are $< R$, and $R<N$.
Therefore, the innermost sum in \eqref{E:Thm2LDef} is $0$
if $(a,[d,e])\ne 1$.
Accordingly, we need an analogue of $\Omega_d(\cH)$ for reduced residue
classes. For squarefree $d$, we define
\begin{equation} \label{E:Omega*Def}
    \Omega_d^*(\cH)= 
     \{a\in \bZ_d: 
	 (a,d)=1 \text{ and } P(a;\cH)\equiv 0 \pmod d \}.
\end{equation}
Let $\nu_d^*(\cH)$ be the cardinality of 
$\Omega_d^*(\cH)$. 
For brevity, we will usually write $\nu_d^*$ in place of 
$\nu_d^*(\cH)$.

When $d_1,d_2$ are squarefree and $(d_1,d_2)=1$, the 
set $\Omega_{d_1d_2}^*(\cH)$ is the image of 
$\Omega_{d_1}^*\times \Omega_{d_2}^*$ under the isomorphism $\xi$
of \eqref{E:CRT}. 
Therefore, the function $\nu^*$ is multiplicative.
Moreover, when $p$ is prime, 
\begin{equation*}
\nu^*_p = \nu_p -1,
\end{equation*}
because we are assuming that $0\in \cH$.

In this context, the most natural analogue of $\fS(\cH)$ is the
product
\begin{equation} \label{E:fS*Def}
    \fS^*(\cH) = \prod_{p} 
	       \left( 1-\frac{\nu^*_p}{p-1} \right)
	       \left(1 - \frac{1}{p} \right)^{-k+1}.
\end{equation}
Note, however that
\begin{align} \label{E:SingularSeriesEq}
       \fS^*(\cH) = & \prod_p 
		      \left( 1 - \frac{\nu_p-1}{p-1} \right)
		      \left( 1- \frac{1}{p} \right)^{-k+1} \\
		  = & \prod_p \left( \frac{p-\nu_p}{p-1} \right)   
			    \left( \frac{p-1}{p} \right)
			    \left( 1- \frac{1}{p} \right)^{-k} \notag \\
		   = & \fS(\cH). \notag
\end{align}

Returning to $\cL$, we write this sum as 
\begin{equation} \label{E:LEst}
    \cL  = \sum_{d,e} \lambda_{d,\ell_1}\lambda_{e,\ell_2} 
              \sum_{a\in \Omega^*_{[d,e]} (\cH)} 
               \sum_{\substack{ N < p \le 2N \\ p\equiv a \pmod {[d,e]}}}
                \log p 
	     = NS+O(T),
\end{equation}	   
where
\begin{equation} \label{E:SDef-Thm6}
S =  \sum_{d,e} 
  \frac{\lambda_{d,\ell_1}\lambda_{e,\ell_2} \nu^*_{[d,e]}}{\phi([d,e])}
\end{equation}
and
\begin{equation*}
T= \sum_{d,e} |\lambda_{d,\ell_1}\lambda_{e,\ell_2}|
		          \nu^*_{[d,e]} E^*(N,[d,e]).
\end{equation*}
By Lemma \ref{L:lambdaBound} and Lemma \ref{L:HR35},
\begin{equation} \label{E:TEst}
    T   \ll (\log R)^{2k+\ell_1+\ell_2} 
             \sumflat_{r< R^2} (3k-3)^{\omega(r)} E^*(N,r)  
        \ll (N/\log N).
\end{equation}

We now consider the sum $S$.
We shall define
\begin{equation} \label{E:Thm2-fDef}
f^*(r) = \frac{\phi(r)}{\nu^*_r}.
\end{equation}
However, we need to take some care with this definition because 
there may be terms with $\nu^*_r=0$. 
However, $\nu^*_p=k-1$ for all but finitely many primes $p$,
so there are at most finitely many primes $p$ such that
$\nu^*_p=0$. We define
\begin{equation} \label{E:ADef}
A= A(\cH) = \prod_{\substack{ p \\ \nu^*_p(\cH)=0}} p,
\end{equation}
and we use the definition in \eqref{E:Thm2-fDef} for any
$r$ with $(r,A)=1$. 
We define $f_1^*$, a function analogous to $f_1$,
by taking
\begin{equation*}
    f_1^*(r)= f^**\mu(r)
\end{equation*}
for $r$ with $(r,A)=1$. 
For future reference, we note that if $p$ is a prime and $p\nmid A$, then
\begin{equation*}
    f^*(p)= \frac{p-1}{\nu_p-1},\,\,\, f^*_1(p) =  \frac{p-\nu_p}{\nu_p-1}. 
\end{equation*}

With this definition of $f^*$, we now have
\begin{equation*}
    S=  \sumprime_{d,e} \frac{\lambda_{d,\ell_1}\lambda_{e,\ell_2}}{f^*(d) f^*(e)}
          \sum_{\substack{r|d\\r|e}} f_1^*(r). 
\end{equation*}
Here, and in the sequel, we use $\sumprime_{\null}$ to denote that the 
sum is over values of the indices that are relatively prime to $A$. 
Interchanging the order of summation, we get
\begin{align}
     S = & \sumprime_r {f_1^*(r)} 
           \left(\sumprime_{d} \frac{\lambda_{dr,\ell_1}}{f^*(dr)} \right) 
         \left(\sumprime_{e} \frac{\lambda_{er,\ell_2}}{f^*(er)} \right) 
	 \label{E:SEst1-Thm6} \\
        = & \sumprime_r \frac{y^*_{r,\ell_1} 
        y^*_{r,\ell_2}}{f_1^*(r)}, \notag
\end{align}
where the
quantity $y_{r,\ell}^*$ is analogous to $y_{r,\ell}$ and is  defined as
\begin{equation} \label{E:y*Def}
    y^*_{r,\ell}=
      \begin{cases}
	\mu(r) f_1^*(r) 
	\displaystyle{
	 \sumprime_d \frac{\lambda_{dr,\ell}}{f^*(dr)} } &
	 \text{ if $(r,A)=1$ and $r<R$,}\\
	 0 &\text{ otherwise}.
      \end{cases}
\end{equation}

Upon using \eqref{E:lambda-yRelation}, 
the original definition of $\lambda_{d,\ell}$, 
we see that
\begin{align*}
    \frac{\mu(r) y_{r,\ell}^*}{f_1^*(r)}  & =
     \sumprime_d \frac{\lambda_{dr,\ell}}{f^*(dr)} =
     \sumprime_d \frac{\mu(dr)}{f^*(dr)} f(dr) 
	\sum_t \frac{y_{rdt,\ell}}{f_1(rdt)} \\
     & = \frac{\mu(r)f(r)}{f^*(r)f_1(r)} 
	 \sumprime_{\substack{d \\ (d,r)=1}} 
	     \frac{\mu(d)f(d)}{f^*(d)}
	  \sum_t \frac{y_{rdt,\ell}}{f_1(dt)}\\
     & = \frac{\mu(r)f(r)}{f^*(r) f_1(r)} 
	    \sum_{\substack{m\\(m,r)=1}}
	     \frac{y_{rm,\ell}}{f_1(m)} 
	     \sumprime_{d|m}\frac{\mu(d)f(d)}{f^*(d)}.
\end{align*}

Note that $m$ can be any squarefree integer; 
we need not have $(m,A)=1$. Now
\begin{align*}
   \sumprime_{d|m}\frac{\mu(d)f(d)}{f^*(d)} = & 
   \prod_{p|m, p\nmid A} \left(1 - \frac{f(p)}{f^*(p)} \right) 
      \\
  = & \prod_{p|m, p\nmid A} \left(\frac{p-\nu_p}{\nu_p(p-1)} \right) \\
  = & \prod_{p|m} \left(\frac{p-\nu_p}{\nu_p(p-1)}\right).
\end{align*}
We may drop the condition that $p\nmid A$ in the last line 
because when $p|A$, $\nu_p=1$, and $ (p-\nu_p)/(\nu_p(p-1))=1$.
Therefore
\begin{equation} \label{E:y*-Formula1}
    \frac{1}{f_1(m)}\sumprime_{d|m}\frac{\mu(d)f(d)}{f^*(d)}=
      \prod_{p|m} \frac{p-\nu_p}{\nu_p(p-1)f_1(p)} =
    \frac{1}{\phi(m)}.
\end{equation}
Moreover,
\begin{equation} \label{E:y*-Formula2}
    \frac{f_1^*(r) f(r)}{f^*(r)f_1(r)} = \frac{r}{\phi(r)}
\end{equation}
when $(r,A)=1$, and so
\begin{equation} \label{E:y*1}
    y^*_{r,\ell} =
    \mu^2(r)
    \frac{\fS(\cH)}{\ell!} \frac{r}{\phi(r)}
	\sum_{\substack{m < R/r \\ (m,r)=1}}
       \frac{\mu^2(m)}{\phi(m)} (\log R/rm)^{\ell}
\end{equation}
when $(r,A)=1$.

For the inner sum, we use Lemma \ref{L:ExtendedWirsing} with
\begin{equation*} 
        \gamma(p)= 
	\begin{cases}
	    1 & \text{ if $p\nmid r$,}\\
	    0 & \text{ if $p|r$}.\\
	 \end{cases}
\end{equation*}
The hypotheses \eqref{E:Omega1} and \eqref{E:Omega2} are satisfied with 
$\kappa=1$, some absolute constants $A_1,A_2$, and 
\begin{equation*}  
    L= \sum_{p|r} \frac{\log p}{p} + O(1).
\end{equation*}
Let 
\begin{equation} \label{E:rhoDef}
    \rho(r)= 1+ \sum_{p|r} \frac{\log p}{p},
\end{equation}
so that $L\ll \rho(r)$.
With this choice of $\gamma$, we have
\begin{equation*}
    c_{\gamma} = \prod_{p|r} \left(1-\frac{1}{p} \right) = 
    \frac{\phi(r)}{r}.
\end{equation*}
We therefore conclude that
\begin{equation} \label{E:PhiSum}
    \sum_{\substack{m< R/r\\ (m,r)=1}} \frac{\mu^2(m)}{\phi(m)} 
    (\log R/rm)^{\ell} =
    \frac{\phi(r)}{r} \frac{(\log R/r)^{\ell+1}}{\ell+1} 
       +O\left( \rho(r) (\log 2R/r)^{\ell} \right).
\end{equation}
We remark parenthetically that Hildebrand \cite{Hildebrand} gave a 
more precise formula for this sum in the case $\ell=0$. It is possible
to use his result to derive a more accurate version of 
\eqref{E:PhiSum}, but the above version is sufficient for our purposes.

From \eqref{E:PhiSum} and \eqref{E:y*1}, we deduce that 
when $(r,A)=1$ and $r<R$,
 \begin{equation} \label{E:y*2}
    y^*_{r,\ell} = 
      \mu^2(r)
      \frac{\fS(\cH)}{(\ell+1)!} (\log R/r)^{\ell+1} 
      + O\left( \frac{\mu^2(r)\rho(r) r}{\phi(r)} \fS(\cH) 
      (\log 2R/r)^{\ell} \right).
\end{equation}
 
We plug this back into our formula for $S$ 
in \eqref{E:SEst1-Thm6} to get 
\begin{equation} \label{E:SDecomp}
    S   = \sumprime_{r < R} 
            \frac{ y^*_{r,\ell_1} y^*_{r,\ell_2} } {f_1^*(r)} 
	  =  V +
	     O\left(\fS(\cH)^2 (\log R)^{\ell_1+\ell_2+1} W\right),
\end{equation}
where 
\begin{equation} \label{E:VDef}
    V=  \frac{\fS(\cH)^2 }{(\ell_1+1)!(\ell_2+1)!}
         \sumprime_{r < R}
          \frac{\mu^2(r)}{f_1^*(r)} 
	   (\log R/r)^{\ell_1+\ell_2+2} 
\end{equation}
and 
\begin{equation} \label{E:WDef}
    W =  
          \sumprime_{r < R}
          \frac{\mu^2(r)}{f_1^*(r)} \frac{\rho(r) r}{\phi(r)}.
\end{equation}         
	  
We will use Lemma \ref{L:ExtendedWirsing} for $V$. 
We will need to estimate a similar sum in Section \ref{S:Thm7},
so it is convenient to have the following lemma that is general enough
to cover both situations.

\begin{lemma} \label{L:ExtendedWirsing-f*}
If $d$ is squarefree, $d < R$, and $a$ is
    a non-negative integer, then
\begin{align*}
   \sumprime_{\substack{r< R/d \\ (r,d)=1}}  
     \frac{\mu^2(r)}{f_1^*(r)} (\log R/dr)^a
	= & \frac{1}{\fS(\cH)} 
	   \frac{a!}{(k+a-1)!} (\log R/d)^{k+a-1} 
	    \prod_{p|d} \left( \frac{p-\nu_p}{p-1} \right)
	   \\
	   & + O\left(
	           (\beta(\cH) + \rho(d)) (\log 2R/d)^{k+a-2}\right).
\end{align*}
\end{lemma}
 
\begin{proof}
We apply Lemma \ref{L:ExtendedWirsing} with     	    
\begin{equation*}
\gamma(p)=
  \begin{cases}
       \displaystyle \frac{p\nu^*_p}{p-1} & \text{if $(p,d)=1$},\\
        0 & \text{if $p|d$}.
   \end{cases}
\end{equation*}
With this definition for $\gamma$, we have
\begin{equation*}
g(p)=\frac{\gamma(p)}{p-\gamma(p)} = \frac{1}{f^*_1(p)}
\end{equation*}	
when $(p,Ad)=1$. Moreover,
\begin{equation*}
\nu^*_p = \nu_p -1 \le \min(k-1,p-2),
\end{equation*}
so \eqref{E:Omega1} is true with $A_1=k$. 
For \eqref{E:Omega2}, we first note that
\begin{align*}
\sum_{w\le p <z} & \frac{\gamma(p) \log p}{p}   = 
\sum_{\substack {w\le p< z \\ (p,d)=1}} \frac{(\nu_p-1)\log p}{p-1} \\
& = (k-1) \sum_{w \le p < z} \frac{\log p}{p-1} -
  \sum_{\substack{w \le p < z \\ (p,d)=1}} 
     \frac{(k-\nu_p)\log p}{p-1} -
  \sum_{\substack{w \le p < z \\ p|d}} 
     \frac{(k-1)\log p}{p-1}.
\end{align*} 
Now
\begin{align*}
\sum_{w\le p < z} & \frac{\log p}{p-1}   =\log(z/w) +O(1),\\
\sum_{\substack{w \le p < z \\ (p,d)=1}} & \frac{(k-\nu_p)\log p}{p-1} 
         \le \beta(\cH)+O(1), \\
\sum_{\substack{w \le p < z \\ p|d}} &
     \frac{(k-1)\log p}{p-1}   \le (k-1) \rho(d),	
\end{align*}
so \eqref{E:Omega2} is satisfied with 
$\kappa=k-1$, $A_2$ some constant depending only on $k$,
and $L= \beta(\cH)+(k-1)\rho(d)+O(1) \ll \beta(\cH)+\rho(d)$.
Finally, we note that in this situation,
\begin{align*}
    c_\gamma & = \prod_{p} \left( 1-\frac{\nu^*_p}{p-1} \right)^{-1}
                      \left( 1-\frac{1}{p} \right)^{k-1}
	      \prod_{p|d} \left( 1- \frac{\nu^*_p}{p-1} \right) \\
	   & = \frac{1}{\fS(\cH)}\prod_{p|d} \left( \frac{p-\nu_p}{p-1} \right)
\end{align*}
by \eqref{E:SingularSeriesEq}.
\end{proof}

From the previous lemma, with $d=1$, we see that 
\begin{align} \label{E:VEst}
V   & = 
          \frac{\fS(\cH)}{(\ell_1+1)!(\ell_2+1)!} 
	  \frac{(\ell_1 + \ell_2 + 2)!}{(k+\ell_1+\ell_2+1)!}
	  (\log R)^{k+\ell_1+\ell_2+1}  \\
      & \phantom{= \fS(\cH)(\ell+1)} 
          +O(\beta(\cH)\fS(\cH)^2 (\log R)^{k+\ell_1+\ell_2}) \notag \\
      & =
         \binomial{\ell_1+\ell_2+2}{\ell_1+1} \fS(\cH)
	  \frac{(\log R)^{k+\ell_1+\ell_2+1}}{(k+\ell_1+\ell_2+1)!} 
	  + O(\beta(\cH)\fS(\cH)^2 (\log R)^{k+\ell_1+\ell_2}). \notag
\end{align}

The sum $W$ may be estimated by relatively trivial means. 
Now
\begin{align} \label{E:WEst}
    W &= \sumprime_{r < R} 
           \frac{\mu^2(r) r}{f_1^*(r)\phi(r)} 
	      \left( 1+ \sum_{p|r} \frac{\log p}{p} \right) 
          \\
        &= \sumprime_{r < R} \frac{\mu^2(r) r}{f_1^*(r)\phi(r)} 
	    + \sumprime_{p < R} \frac{\log p}{p} 
	      \sumprime_{\substack{r<R \\ p|r}}
	              \frac{\mu^2(r) r}{f_1^*(r)\phi(r)} \notag
	   \\
 	 &=\sumprime_{r < R} \frac{\mu^2(r) r}{f_1^*(r)\phi(r)} 
	    + \sumprime_{p < R} \frac{\log p}{f^*_1(p)\phi(p)} 
 	      \sumprime_{\substack{r<R/p\\(r,p)=1}}
 	              \frac{\mu^2(r) r}{f_1^*(r)\phi(r)} \notag
 	    \\
 	 & \ll 
 	    \left( 1+ \sumprime_{p<R} \frac{\log p}{f^*_1(p)\phi(p)} \right)
 	       W^* \ll W^*, \notag
\end{align}
where 
\begin{equation} \label{E:W*Def}
    W^*= \sumprime_{r < R} 
           \frac{\mu^2(r) r}{f_1^*(r)\phi(r)}
	 =\sumflat_{r < R} \frac{\nu^*_r h(r)}{r},
\end{equation}
and 
\begin{equation} \label{E:hDef}
    h(r) = \prod_{p|r} \frac{p^2}{(p-\nu_p)(p-1)}.
\end{equation}
Let $h_1=h*\mu$, so that
\begin{equation*}
    h_1(d)=
      \prod_{p|d} \frac{p(\nu_p+1)-\nu_p}{(p-1)(p-\nu_p)}.
\end{equation*}
Then
\begin{equation} \label{E:W*EstA}
    W^*  = \sumflat_{r < R} \frac{\nu^*_r}{r} \sum_{d|r} h_1(d)  
          = \sumflat_{d < R} \frac{h_1(d) \nu^*_d}{d}
	     \sumflat_{\substack{r < R/d \\ (r,d)=1}} \frac{\nu^*_r}{r} 
	\le  \prod_{p< R} \left(1+ \frac{h_1(p) \nu^*_p}{p} \right)
	     \sumflat_{r < R} \frac{\nu^*_r}{r}.
\end{equation}    
The sum on the right-hand side of \eqref{E:W*EstA} is $\ll (\log R)^{k-1}$
by Lemma \ref{L:HR34}. The product is $\ll 1$ because 
\begin{equation*}
    \sum_{p< R} \log \left(1+\frac{h_1(p)\nu^*_p}{p}\right) \ll
    \sum_{p< R} \frac{\nu_p^2}{(p-1)(p-\nu_p)} \ll 1.
\end{equation*}
We conclude that  $ W^*\ll (\log R)^{k-1}$, and so
\begin{equation} \label{E:WEstFinal}
    W \ll   (\log R)^{k-1}.
\end{equation}
Combining the above with the estimate in \eqref{E:VEst}
gives
\begin{equation} \label{E:SEst}
    S=   \binomial{\ell_1+\ell_2+2}{\ell_1+1} \fS(\cH)
	  \frac{(\log R)^{k+\ell_1+\ell_2+1}}{(k+\ell_1+\ell_2+1)!} 
	  + O\left(
	       \beta(\cH)\fS(\cH)^2 (\log R)^{k+\ell_1+\ell_2}
	         \right).
\end{equation}

The first part of Theorem 6 (statement \eqref{E:Thm6Part1}) 
now follows by combining
\eqref{E:LEst},
\eqref{E:TEst}, and
\eqref{E:SEst}.

\section{Proof of Theorem \ref{T:Thm6}, Part 2}\label{S:Thm6-P2}
In this section, we consider Theorem \ref{T:Thm6} in the case
$h_0\notin \cH$.
As in the previous section, our problem is translation invariant, 
so we we may assume that $0\notin \cH$ and $\cH^{0}=\cH\union\{0\}$.
Consequently, $P(n;\cH^{0})=nP(n;\cH)$.

Now let $\cL$ be the left-hand side of \eqref{E:Thm6Part2}.
If $n$ is a prime with $ N< n \le 2N$, then $1$ is the only divisor of 
$n$ less than $N$. When $d< R<N$, we have
$d|P(n;\cH)$ if and only if $d|P(n;\cH^0)$. 
Consequently,
\begin{equation} \label{E:Thm2P2LDef}
    \cL=\sum_{d,e} \lambda_{d,\ell_1} \lambda_{e,\ell_2} 
     \sum_{\substack{N < n \le 2N\\ [d,e]|P(n;\cH^{0})}} \varpi(n)  =
    \sum_{d,e} \lambda_{d,\ell_1} \lambda_{e,\ell_2}
      \sum_{a\in \Omega^*_{[d,e]}(\cH^{0})} 
      \sum_{\substack{N < p \le 2N \\ p\equiv a \pmod{[d,e]} }}
      \log p.
\end{equation}

In parallel to the argument in \eqref{E:LEst} through 
\eqref{E:TEst}, we find that 
\begin{equation*}
    \cL= N S + T,
\end{equation*}
where
\begin{equation}\label{E:SDef-Thm6P2}
    S =  \sum_{d,e} 
      \frac{\lambda_{d,\ell_1}\lambda_{e,\ell_2} \nu^*_{[d,e]}(\cH^0)}{\phi([d,e])}
\end{equation}
and
\begin{equation*}
    T= \sum_{d,e} |\lambda_{d,\ell_1}\lambda_{e,\ell_2}|
			      \nu^*_{[d,e]}(\cH^0) E^*(N,[d,e])
    \ll N/\log N.
\end{equation*}
Therefore
\begin{equation} \label{E:Thm2P2-Step1}
    \cL= NS +O(N/\log N). 
\end{equation}
The rest of this section is devoted to evaluating the sum $S$.

For brevity, we write $\nu_r^\dagger$ for 
$\nu_r^*(\cH^0)$. 
Let 
\begin{equation*}
    A_0= A(\cH^0)= \prod_{\substack{p\\ \nu_p^\dagger=0}} p.
\end{equation*}
For squarefree $r$ with $(r,A_0)=1$, we define
\begin{equation} \label{E:f+Def}
    f^\dagger(r) = \frac{\phi(r)}{\nu_r^\dagger} = 
     \prod_{p|r} \left(\frac{p-1}{\nu_p^\dagger} \right),
\end{equation}
and 
\begin{equation} \label{E:f1+Def}
    f_1^\dagger(r) =f^\dagger*\mu(r) = 
    \prod_{p|r} \left(\frac{p-1-\nu_p^\dagger}{\nu_p^\dagger} \right)
\end{equation}

Note that
\begin{equation*}
    \nu_p^\dagger =
    \begin{cases}
	 \nu_p & \text{if } 0\notin \Omega_p(\cH) \\
	 \nu_p-1 & \text{if } 0\in  \Omega_p(\cH) .
     \end{cases}
\end{equation*}
We are assuming that $0\notin\cH$, so there are only finitely many 
primes $p$ with $0\in \Omega_p(\cH)$. 
Let
\begin{equation} \label{E:B0Def}
    B_0= B_0(\cH)=\prod_{\substack{p\\ \nu_p^\dagger=\nu_p-1}} p
		 =\prod_{\substack{p\\ 0\in \Omega_p(\cH)}} p.
\end{equation}
In fact, $0\in \Omega_p(\cH)$ if and only if $p$ divides $h$
for some $h\in \cH$. Therefore $B_0$ is the squarefree kernel
of the product of all elements of $\cH$.

For future reference, we note that when $(r,A_0)=1$,
\begin{equation*}
    f^\dagger(r)=
    \prod_{\substack{p|r\\p\nmid B_0}} \left(\frac{p-1}{\nu_p}\right)
    \prod_{\substack{p|r\\ p|B_0}} \left(\frac{p-1}{\nu_p-1}\right)
\end{equation*}
and 
\begin{equation*}
    f_1^\dagger(r)=
	\prod_{\substack{p|r\\p\nmid B_0}} 
	\left(\frac{p-1-\nu_p}{\nu_p}\right)
	\prod_{\substack{p|r\\ p|B_0}} 
	 \left(\frac{p-\nu_p}{\nu_p-1}\right).
\end{equation*}

With the above definitions of $f^\dagger$ and $f_1^\dagger$, we may 
write
\begin{align*}
    S = &   \sumprime_{d,e} 
      \frac{\lambda_{d,\ell_1}\lambda_{e,\ell_2} }{f^\dagger([d,e])} \\
      = & \sumprime_{d,e} 
      \frac{\lambda_{d,\ell_1}\lambda_{e,\ell_2}}{f^\dagger(d)f^\dagger(e)}
      \sum_{\substack{r|d\\r|e}} f_1^\dagger(r) \\
      = & \sumprime_{r} f_1^\dagger(r)
	\left(\sumprime_d\frac{\lambda_{dr,\ell_1}}{f^\dagger(dr)}\right)
	\left(\sumprime_e\frac{\lambda_{er,\ell_2}}{f^\dagger(er)}\right),
\end{align*}
where $\sumprime_\null$ denotes that the sum is over values of the 
indices that are relatively prime to $A_0$.
We get
\begin{equation}\label{E:SEst1-Thm6P2}
    S=\sumprime_r 
    \frac{y^\dagger_{r,\ell_1}y^\dagger_{r,\ell_2}}
	 {f_1^\dagger(r)},
\end{equation}
where we define
\begin{equation} \label{E:ydaggerDef}
    y^\dagger_{r,\ell} = 
    \begin{cases}
    \mu(r) f_1^\dagger(r) 
    \displaystyle{
     \sumprime_d \frac{\lambda_{dr,\ell}}{f^\dagger(dr)} } &
     \text{ if $(r,A_0)=1$ and $r<R$,}\\
     0 & \text{ otherwise.}
     \end{cases}
\end{equation}

Upon using \eqref{E:lambda-yRelation}, 
our original definition of $\lambda_{d,\ell}$, 
we see that
\begin{align*}
    \frac{\mu(r) y_{r,\ell}^\dagger}{f_1^\dagger(r)}  & =
     \sumprime_d \frac{\lambda_{dr,\ell}}{f^\dagger(dr)} =
     \sumprime_d \frac{\mu(dr)}{f^\dagger(dr)} f(dr) 
	\sum_t \frac{y_{rdt,\ell}}{f_1(rdt)} \\
     & = \frac{\mu(r)f(r)}{f^\dagger(r)f_1(r)} 
	 \sumprime_{\substack{d \\ (d,r)=1}} 
	     \frac{\mu(d)f(d)}{f^\dagger(d)}
	  \sum_t \frac{y_{rdt,\ell}}{f_1(dt)}\\
     & = \frac{\mu(r)f(r)}{f^\dagger(r) f_1(r)} 
	    \sum_{\substack{m\\(m,r)=1}}
	     \frac{y_{rm,\ell}}{f_1(m)} 
	     \sumprime_{d|m}\frac{\mu(d)f(d)}{f^\dagger(d)}.
\end{align*}

Now 
\begin{equation*}
    \sumprime_{d|m}\frac{\mu(d)f(d)}{f^\dagger(d)}=
    \prod_{\substack{p|m\\p\nmid A_0}} 
      \left( 1- \frac{f(p)}{f^\dagger(p)} \right)
     =
     \prod_{p|m}
     \left( 1- \frac{ p\nu_p^\dagger}{(p-1)\nu_p} \right).
\end{equation*}
The condition $p\nmid A_0$ can be dropped because 
$\nu^\dagger_p=0$ when $p|A_0$. 
Therefore
\begin{align*}
    \sumprime_{d|m}\frac{\mu(d)f(d)}{f^\dagger(d)}
    = &
    \prod_{\substack{p|m\\p\nmid B_0}} 
       \left( 1- \frac{p}{p-1} \right)
     \prod_{\substack{p|m\\ p|B_0}}
       \left( 1 - \frac{p(\nu_p-1)}{\nu_p(p-1)} \right) \\
     = & \frac{\mu(m)}{\phi(m)} f_2(m),
\end{align*}
where $f_2$ is the multiplicative function defined by
\begin{equation} \label{E:f2Def}
    f_2(p) =
    \begin{cases}
	1 & \text{ if } p\nmid B_0,\\
	-f_1(p) & \text{ if } p|B_0.
     \end{cases}
\end{equation}
In other words,
\begin{equation*}
    f_2(m) = \mu((m,B_0)) f_1((m,B_0)).
\end{equation*}

Therefore
\begin{align} 
    y^\dagger_{r,\ell} = &
    \mu^2(r)
    \frac{f_1^\dagger(r) f(r)}{f^\dagger(r) f_1(r)}
    \sum_{\substack{m< R/r\\(m,r)=1}} 
    \frac{y_{rm,\ell}}{f_1(m)} \frac{\mu(m)}{\phi(m)} f_2(m) 
    \label{E:y+Step1}
    \\
    = & 
    \mu^2(r)
    \frac{\fS(\cH)}{\ell!} 
    \frac{f_1^\dagger(r) f(r)}{f^\dagger(r) f_1(r)}
       \sum_{\substack{m< R/r\\(m,r)=1}}
       \frac{\mu(m) f_2(m)(\log R/rm)^\ell}{f_1(m)\phi(m)}.
       \notag
\end{align}

The sum
\begin{equation*}
     \sum_{\substack{m=1\\(m,r)=1}}^\infty
       \frac{\mu(m) f_2(m)}{f_1(m)\phi(m)}
\end{equation*}
converges, and so one would expect that
\begin{equation*}
    y^\dagger_{r,\ell} \asymptotic
    \mu^2(r) 
    \frac{\fS(\cH)}{\ell!} (\log R/r)^\ell
	\frac{f_1^\dagger(r) f(r)}{f^\dagger(r) f_1(r)}
	   \sum_{\substack{m=1\\(m,r)=1}}^{\infty}
	   \frac{\mu(m) f_2(m)}{f_1(m)\phi(m)}
\end{equation*}
when $r<R$ and $(r,A_0)=1$.
From Lemma \ref{L:y+InfSum} below, we would then obtain
\begin{equation*} % \label{E:y+Conjecture}
    y^\dagger_{r,\ell} \asymptotic
       \mu^2(r) \frac{\fS(\cH^0)}{\ell!}
       (\log R/r)^\ell,
\end{equation*}
and we will ultimately prove this.
This asymptotic relation 
should be compared to
\eqref{E:yrChoice} and \eqref{E:y*2}.

\begin{lemma} \label{L:y+InfSum}
    If $r$ is squarefree and $(r,A_0)=1$, then
\begin{equation*}
    \frac{f_1^\dagger(r) f(r)}{f^\dagger(r) f_1(r)}
	  \sum_{\substack{m=1\\(m,r)=1}}^\infty
	  \frac{\mu(m) f_2(m)}{f_1(m)\phi(m)}=
     \frac{\fS(\cH^0)}{\fS(\cH)}.
\end{equation*}
\end{lemma}

\begin{proof} For $r$ satisfying our hypotheses,
it is convenient to define
\begin{equation}\label{E:FGDefs}
    F(r)=  \frac{f_1^\dagger(r) f(r)}{f^\dagger(r) f_1(r)}
    \text{ and }
    G(r)=   \sum_{\substack{m=1\\(m,r)=1}}^\infty
	  \frac{\mu(m) f_2(m)}{f_1(m)\phi(m)},
\end{equation}
so that the left-hand side of the proposed result is $F(r)G(r)$.
We begin by noting that
\begin{equation*}
   F(r)=\prod_{p|r} F(p)=\prod_{p|r} \frac{p(p-1-\nu_p^\dagger)}{(p-1)(p-\nu_p)}.
\end{equation*}
Moreover,
\begin{align*}
   G(r)
    = & \prod_{p\nmid r}
	 \left(1-\frac{f_2(p)}{\phi(p)f_1(p)}\right) \\
    = & \prod_{\substack{p\nmid B_0\\p\nmid r}} 
	\frac{p(p-1-\nu_p)}{(p-1)(p-\nu_p)}
       \prod_{\substack{p|B_0\\p\nmid r}}
	 \frac{p}{p-1}\\
    = & \prod_{p\nmid r} 
	\frac{p(p-1-\nu_p^\dagger)}{(p-1)(p-\nu_p)}
    =  \prod_{p\nmid r} F(p).
\end{align*}
In the last line, we used the fact that $\nu_p^\dagger=\nu_p$ if 
$p\nmid B_0$ and $\nu_p^\dagger=\nu_p-1$ if $p\mid B_0$.
Combining the last two results yields
\begin{equation}\label{E:y+InfSum-Step1}
F(r)G(r)=\prod_{p}  \frac{p(p-1-\nu_p^\dagger)}{(p-1)(p-\nu_p)}
=\prod_p F(p).
\end{equation}

On the other hand, if we replace $\cH$ by $\cH^0$ and $k$ by $k+1$ 
in \eqref{E:SingularSeriesEq}, then we obtain
\begin{equation*}
\fS(\cH^0)=\fS^*(\cH^0)=
\prod_p  \left(1-\frac{\nu_p^\dagger}{p-1}\right)
\left(1-\frac{1}{p}\right)^{-k}.
\end{equation*}
We combine this with the definition of $\fS(\cH)$ given
in \eqref{E:fSDefinition} to get
\begin{equation} \label{E:F=S}
    \frac{\fS(\cH^0)}{\fS(\cH)}=
    \prod_p \frac{p(p-1-\nu_p^\dagger)}{(p-1)(p-\nu_p)}
    =\prod_p F(p).
\end{equation}
The lemma follows by comparing this with \eqref{E:y+InfSum-Step1}.
\end{proof}

\begin{lemma} \label{L:y+2}
Suppose $\ell\ge 1$.
If  $r< R$ and  $(r,A_0)=1$, then
\begin{equation} \label{E:y+2}
    y^\dagger_{r,\ell} = 
     \mu^2(r)
      \frac{\fS(\cH^0)}{\ell!} (\log  R/r)^\ell
      + O\left(\mu^2(r)
	   \beta(\cH^0)\fS(\cH^0)(\log 2R/r)^{\ell-1}
	   \right). 
\end{equation}
\end{lemma}

\begin{proof} 
From the definition of $y^\dagger_{r,\ell}$ in \eqref{E:ydaggerDef},
the lemma is trivial if $r$ is not squarefree.
For the remainder of the proof, we assume that $r$ is squarefree,
$(r,A_0)=1$, and $r< R$.

We start from the expression for 
$y^\dagger_{r,\ell}$ given in \eqref{E:y+Step1}.
For a given $m$ in the inner sum, write $m=\delta n$, where
$\delta|B_0$ and $(n,B_0)=1$. Then $f_2(m)=\mu(\delta)f_1(\delta)$
and 
\begin{equation*}
    \frac{\mu(m)f_2(m)}{f_1(m)\phi(m)}=
    \frac{\mu^2(\delta)\mu(n)}{\phi(\delta)\phi(n)f_1(n)}.
\end{equation*}
Therefore \eqref{E:y+Step1} may be transformed into
\begin{equation*} 
    y^\dagger_{r,\ell}=
    \frac{\fS(\cH)F(r)}{\ell!}
    \sum_{\substack{\delta|B_0\\(\delta,r)=1}}
      \frac{\mu^2(\delta)}{\phi(\delta)}
     \sum_{\substack{n<R/r\delta \\ (n,rB_0)=1}}
       \frac{\mu(n)}{f_1(n)\phi(n)} (\log R/r\delta n)^\ell.
\end{equation*}
If we  set 
\begin{equation}\label{E:B1Def}
    B_1= \prod_{\substack{p|B_0\\p\nmid r}} p = \frac{B_0}{(B_0,r)},
\end{equation}
then the above equation for $y^\dagger_{r,\ell}$ may be written as
\begin{equation}\label{E:y+Step2}
    y^\dagger_{r,\ell}=
       \frac{\fS(\cH)F(r)}{\ell!}
       \sum_{\delta|B_1}
	 \frac{\mu^2(\delta)}{\phi(\delta)}
	\sum_{\substack{n<R/r\delta \\ (n,rB_1)=1}}
	  \frac{\mu(n)}{f_1(n)\phi(n)} (\log R/r\delta n)^\ell. 
\end{equation}
For future reference, note that $B_0|rB_1$.

Now let 
\begin{equation}\label{E:YDef}
    Y(x;d,\ell)=
     \sum_{\substack{n<x\\(n,d)=1}} 
      \frac{\mu(n)}{f_1(n)\phi(n)} (\log x/n)^\ell,
\end{equation}
so that the innermost sum in \eqref{E:y+Step2} is 
$Y(R/r\delta ; rB_1, \ell)$.

Now assume that $\ell\ge 1$. 
We begin our analysis of $Y$ by writing
\begin{align} 
    Y(x;d,\ell)= &
     \sum_{\substack{n< x\\ (n,d)=1}} 
      \frac{\mu(n)}{f_1(n)\phi(n)} 
	\int_n^x \ell (\log x/u)^{\ell-1} \frac{du}{u}
	\label{E:YStep1} \\
     = & \int_1^x \frac{\ell (\log x/u)^{\ell-1}}{u}
	    \sum_{\substack{n< u\\ (n,d)=1}}
	     \frac{\mu(n)}{f_1(n)\phi(n)} du
	     \notag \\
     = & Y_1(x;d,\ell)-Y_2(x;d,\ell),
	     \notag
\end{align}
where
\begin{equation} \label{E:Y1Def}
    Y_1(x;d,\ell)= 
    \int_1^x \frac{\ell(\log x/u)^{\ell-1}}{u}
		\sum_{\substack{n=1\\ (n,d)=1}}^\infty
		 \frac{\mu(n)}{f_1(n)\phi(n)} du,
\end{equation}
and 
\begin{equation} \label{E:Y2Def}
    Y_2(x;d,\ell)= 
    \int_1^x  \frac{\ell (\log x/u)^{\ell-1}}{u}
		\sum_{\substack{n\ge u \\ (n,d)=1}}
		 \frac{\mu(n)}{f_1(n)\phi(n)} du.
\end{equation}

We see immediately that
\begin{equation*}
    Y_1(x;d,\ell)= (\log x)^\ell
    \prod_{p\nmid d} \left( 1- \frac{1}{f_1(p)\phi(p)}\right).
\end{equation*}
If we assume that $B_0|d$, then we may write
\begin{equation}  \label{E:Y1Est}
    Y_1(x;d,\ell)= (\log x)^\ell \prod_{p\nmid d} F(p).
\end{equation}

For $Y_2(x;d,\ell)$ we bound the sum inside the integrand as
\begin{equation} \label{E:Y2Step1}
    \left|
     \sum_{\substack{n\ge u\\ (n,d)=1}}
      \frac{\mu(n)}{f_1(n)\phi(n)}
       \right|
       \le 
       \sum_{n \ge u} 
	\frac{\mu^2(n)}{f_1(n)\phi(n)}
       =\int_u^\infty 
	 \left( 
	  \sum_{u\le n < v} \frac{\mu^2(n)n}{f_1(n)\phi(n)} 
	  \right)
	  \frac{dv}{v^2}.
\end{equation}
Now let 
\begin{equation}\label{E:W+Def}
    W^\dagger(v) = \sum_{ n < v} \frac{\mu^2(n)n}{f_1(n)\phi(n)}.
\end{equation}
This sum is very similar to the sum $W^*$ defined in \eqref{E:W*Def}; 
in fact,
\begin{equation*}
    W^\dagger(v)= \sumflat_{n< v} \frac{\nu_n h(n)}{n},
\end{equation*}
where $h$ was defined in \eqref{E:hDef}. 
We have, similarly to \eqref{E:W*EstA},
\begin{equation*} 
W^\dagger(v) = \sumflat_{n< v} \frac{\nu_n}{n}\sum_{d|n} h_1(d)
 = \sumflat_{d< v} \frac{h_1(d)\nu_d}{d} 
   \sumflat_{\substack{n<v/d \\ (n,d)=1}} \frac{\nu_n}{n}
 \le \prod_{p< n} \left( 1+ \frac{h_1(p)\nu_p}{p}\right)
     \sumflat_{n<v} \frac{\nu_n}{n}.
\end{equation*}
The sum on the right-hand side is $\ll (\log 2v)^k$ by Lemma 
\ref{L:HR34}. 
The product on the right hand side is $\ll 1$ because
\begin{equation*}
    \sum_{p<v} \log\left( 1+ \frac{h_1(p)\nu_p}{p} \right)
    \ll \sum_{p<v} \frac{\nu_p^2}{(p-1)(p-\nu_p)}\ll 1.
\end{equation*}
Therefore
\begin{equation}\label{E:W+Est}
W^\dagger(v) \ll (\log 2v)^k.
\end{equation}

Now we use \eqref{E:W+Est} in \eqref{E:Y2Step1} to get
\begin{equation*}
    \left|
	 \sum_{\substack{n\ge u\\ (n,d)=1}}
	  \frac{\mu(n)}{f_1(n)\phi(n)}
	   \right|
    \ll \int_u^\infty \frac{(\log 2v)^k}{v^2} dv 
    \ll \frac{(\log 2u)^k}{u}.
\end{equation*}
We use this in \eqref{E:Y2Def} to get
\begin{equation}\label{E:Y2EstF}
    Y_2(x;d,\ell)\ll (\log 2x)^{\ell-1} \int_1^x (\log 2v)^k \frac{dv}{v^2}
	      \ll (\log 2x)^{\ell-1}.
\end{equation}
Combining this with \eqref{E:Y1Est} gives
\begin{equation} \label{E:YEst}
    Y(x;d,\ell)= 
    (\log x)^\ell \prod_{p\nmid d} F(p) + O((\log 2x)^{\ell-1})
\end{equation}
when $B_0|d$.

Now we use \eqref{E:YEst} with $d=rB_1$ in \eqref{E:y+Step2} to obtain
\begin{align}   
    y^\dagger_{r,\ell}= 
    \frac{\fS(\cH)}{\ell!} & \left(\prod_{p\nmid B_1} F(p)\right)
    \sum_{\substack{\delta|B_1\\ \delta< R/r}} 
       \frac{\mu^2(\delta)}{\phi(\delta)} (\log R/r\delta)^\ell
       \label{E:y+Step3} \\
      & +O\left(
	    \fS(\cH) F(r)
	   \sum_{\substack{\delta|B_1\\ \delta< R/r}} 
	   \frac{\mu^2(\delta)}{\phi(\delta)}
	   (\log 2R/r\delta)^{\ell-1}
	   \right).
	   \notag
\end{align}

The error term in \eqref{E:y+Step3} is 
\begin{align*}
    \ll & \fS(\cH) F(r) (\log 2R/r)^{\ell-1} \sum_{\delta|B_1} 
	   \frac{\mu^2(\delta)}{\phi(\delta)} \\
     \ll & \fS(\cH) \left( \prod_{p|rB_1} F(p) \right) 
	    (\log 2R/r)^{\ell-1} \\
     \ll & \fS(\cH^0) (\log 2R/r)^{\ell-1}
	   \left( \prod_{p\nmid rB_1} F(p) \right)^{-1}.
\end{align*}
We have used \eqref{E:F=S} in the last line.
Now when $p\nmid B_0$,
\begin{equation*}
    F(p)^{-1} = \left(1 - \frac{\nu_p}{(p-1)(p-\nu_p)} \right)^{-1}
	      = 1+O(1/p^2),
\end{equation*}
so 
\begin{equation*}
    \left( \prod_{p\nmid rB_1} F(p) \right)^{-1} \ll 1.
\end{equation*}
Therefore the error term in \eqref{E:y+Step3} is
\begin{equation} \label{E:y+ErrorEst}
    \ll \fS(\cH^0) (\log 2R/r)^{\ell-1}.
\end{equation}

Now we consider the main term in \eqref{E:y+Step3}, which we write as
\begin{equation} \label{E:y+MainDecomp}
    \frac{\fS(\cH)}{\ell!} \left(\prod_{p\nmid B_1} F(p) \right)
       \left\{ M_1-M_2-M_3\right\},
\end{equation}
where
\begin{align*}
    M_1 = & (\log R/r)^\ell \sum_{\delta|B_1} 
    \frac{\mu^2(\delta)}{\phi(\delta)}, \\
    M_2 = & (\log R/r)^\ell 
      \sum_{\substack{\delta|B_1 \\ \delta\ge R/r}}
     \frac{\mu^2(\delta)}{\phi(\delta)},\\
    M_3 = &\sum_{\substack{\delta|B_1 \\ \delta < R/r}}
    \frac{\mu^2(\delta)}{\phi(\delta)}
    \left\{ (\log R/r)^\ell - (\log R/r\delta)^\ell \right\}.
\end{align*}

For $M_1$, we note that 
\begin{equation*}
    \sum_{\delta|B_1} \frac{\mu^2(\delta)}{\phi(\delta)}
    = \prod_{p|B_1} \frac{p}{p-1}
    = \prod_{p|B_1} F(p) .
\end{equation*}
Therefore
\begin{equation} \label{E:M1Est}
    \frac{\fS(\cH)}{\ell!} M_1 \prod_{p\nmid B_1} F(p)
    = \frac{\fS(\cH)}{\ell!} (\log R/r)^\ell \prod_{p} F(p)
    = \frac{\fS(\cH^0)}{\ell!} (\log R/r)^\ell.
\end{equation}
by \eqref{E:F=S}.

For $M_2$, we note that 
\begin{equation*}
    \sum_{\substack{\delta|B_1 \\ \delta\ge R/r}}
	 \frac{\mu^2(\delta)}{\phi(\delta)}
	\ll \sum_{\delta|B_1} \frac{\mu^2(\delta)}{\phi(\delta)}
	   \frac{\log \delta}{(\log 2R/r)},
\end{equation*}
and 
\begin{align*}
    \sum_{\delta|B_1} \frac{\mu^2(\delta)}{\phi(\delta)} 
    \log \delta = &
    \sum_{\delta|B_1} \frac{\mu^2(\delta)}{\phi(\delta)} 
     \sum_{p|\delta}\log p =
     \sum_{p|B_1} \frac{\log p}{p-1} 
     \sum_{\delta|B_1/p} \frac{\mu^2(\delta)}{\phi(\delta)}\\
     = & \sum_{p|B_1} \frac{\log p}{p-1} 
	    \frac{B_1/p}{\phi(B_1/p)} 
     = \frac{B_1}{\phi(B_1)} \sum_{p|B_1} \frac{\log p}{p}\\
     = & F(B_1) \sum_{p|B_1}  \frac{\log p}{p}.
\end{align*}
Now if $p|B_1$, then $p|B_0$ and $\nu_p(\cH^0)\le k.$
Therefore
\begin{equation*}
    \sum_{p|B_1}\frac{\log p}{p} \le 
    \sum_p \frac{(k+1-\nu_p(\cH^0))\log p}{p}
    =\beta(\cH^0).
\end{equation*}
Consequently,
\begin{equation}\label{E:LogSum}
    \sum_{\delta|B_1} \frac{\mu^2(\delta)}{\phi(\delta)} 
	\log \delta \ll F(B_1) \beta(\cH^0),
\end{equation}
and so
\begin{align}
    \frac{\fS(\cH)}{\ell!} M_2 \prod_{p\nmid B_1} F(p)
     \ll & (\log 2R/r)^{\ell-1} \fS(\cH)\beta(\cH^0) F(B_1)
	  \prod_{p\nmid B_1} F(p) 
	  \label{E:M2Est} \\
     \ll &  \fS(\cH^0) \beta(\cH^0) (\log 2R/r)^{\ell-1}.
     \notag
\end{align}

For $M_3$, we note that when $\delta\le R/r$,
\begin{align*}
    (\log R/r)^\ell - & (\log R/r\delta)^\ell 
     \\
     = &  (\log \delta) 
     \left\{ (\log R/r)^{\ell-1}+ (\log R/r\delta)(\log R/r)^{\ell-2}
	      + \ldots + (\log R/r\delta)^{\ell-1}
      \right\} \\
      \ll & (\log \delta) (\log R/r)^{\ell-1}.
\end{align*}
Thus 
\begin{equation*}
    M_3 \ll (\log 2R/r)^{\ell-1} 
	   \sum_{\delta|B_1} \frac{\mu^2(\delta)}{\phi(\delta)} 
	      \log \delta
	 \ll (\log 2R/r)^{\ell-1} F(B_1) \beta(\cH^0)
\end{equation*}
by \eqref{E:LogSum}, and so
\begin{equation} \label{E:M3Est}
    \frac{\fS(\cH)}{\ell!} M_3 \prod_{p\nmid B_1} F(p)
    \ll (\log 2R/r)^{\ell-1} \fS(\cH^0) \beta(\cH^0).
\end{equation}

Combining the estimates \eqref{E:y+ErrorEst},\eqref{E:M1Est}, 
\eqref{E:M2Est}, and \eqref{E:M3Est} gives the proof of 
Lemma \ref{L:y+2}.
\end{proof}

In reference to the above lemma, we remark that with a bit more
work we could give an estimate valid for $y_{r,0}$ with a somewhat
weaker error term. However, we omit this because it is not necessary 
for the proof of Theorem \ref{T:GPY1}.

We can now complete the estimate of $S$. 
From \eqref{E:SDef-Thm6P2} and Lemma \ref{L:y+2}, we see that
\begin{equation} \label{E:SEst2-T6P2}
    S =
     V^\dagger
    + O\left(\fS(\cH^0)^2\beta(\cH^0)(\log 
    R)^{\ell_1+\ell_2-1}W^\dagger\right),
\end{equation}
where
\begin{equation*}
    V^\dagger =
     \frac{\fS(\cH^0)^2}{\ell_1!\ell_2!}
     \sumprime_{r< R}
     \frac{\mu^2(r)}{f_1^\dagger(r)} (\log R/r)^{\ell_1+\ell_2},
\end{equation*}
and 
\begin{equation*}
    W^\dagger= \sumprime_{r< R}
     \frac{\mu^2(r)}{f_1^\dagger(r)}. 
\end{equation*}
Now $V^\dagger$ is the same as the sum $V$ in \eqref{E:VDef} except 
that $\cH$ has been replaced by $\cH^0$, $k$ has been replaced by 
$k+1$, and $\ell_1,\ell_2$ have been replaced by $\ell_1-1,\ell_2-1$
respectively. 
From \eqref{E:VEst}, we see that 
\begin{equation} \label{E:VEst+}
V^\dagger  =
	 \binomial{\ell_1+\ell_2}{\ell_1} \fS(\cH^0)
	  \frac{(\log R)^{k+\ell_1+\ell_2}}{(k+\ell_1+\ell_2)!} 
	  + O(\beta(\cH^0)\fS(\cH^0)^2 (\log R)^{k+\ell_1+\ell_2-1}).
\end{equation}  
For $W^\dagger$, we use Lemma \ref{L:ExtendedWirsing-f*}
with $a=0$, $d=1$, $f^*$ replaced by $f^\dagger$, $k$ replaced by $k+1$
to get
\begin{equation} \label{E:WEst+}
W^\dagger \ll (\log R)^k.
\end{equation}

Now we combine \eqref{E:SEst2-T6P2}, \eqref{E:VEst+}, and 
\eqref{E:WEst+} to get 
\begin{equation} \label{E:SEst3-T6P2}
    S= 
    \binomial{\ell_1+\ell_2}{\ell_1} \fS(\cH^0)
	  \frac{(\log R)^{k+\ell_1+\ell_2}}{(k+\ell_1+\ell_2)!} 
	  + O(\beta(\cH^0)\fS(\cH^0)^2 (\log R)^{k+\ell_1+\ell_2-1}).
\end{equation}

Equation \eqref{E:Thm6Part2} now follows by combining this 
with \eqref{E:Thm2P2-Step1}.

\section{Proof of Theorem \ref{T:Thm7}} \label{S:Thm7}
We may again assume, without loss of generality, that $h_0=0$. 
Accordingly, we assume throughout this section that $0\in \cH$.

Let $\cL$ denote the sum on the left-hand side in the statement of 
Theorem \ref{T:Thm7}. Then
\begin{equation}  \label{E:LDef}
\cL= \sum_{d,e} \lambda_{d,\ell_1} \lambda_{e,\ell_2} 
   \sum_{a\in \Omega_{[d,e]}(\cH)} 
   \sum_{\substack{N < n \le 2N \\ n \equiv a \pmod {[d,e]}}}
    \varpi*\varpi(n).
\end{equation}
 
In this sum, we have $d,e< R < \sqrt{N}$, so $[d,e]$ has
no prime divisors exceeding $\sqrt{N}$. On the other hand,
if $N< n \le 2N$ and $\varpi*\varpi(n)>0$, then 
$n$ is a product of two primes, at least one of which must exceed
$\sqrt N$. Therefore, the inner sum in \eqref{E:LDef} will be 0
unless $(a,[d,e])=1$ or $(a,[d,e])=p$ for some prime $p< R$.

We write
\begin{equation*} % \label{E:LDecomp}
    \cL=\cL_1+\cL_2,
\end{equation*}
where $\cL_1$ is the sum in \eqref{E:LDef} with the extra condition
that $(a,[d,e])=1$, and
$\cL_2$ is the sum in \eqref{E:LDef} with the extra condition that
$(a,[d,e])=p$ for some prime $p$.

Before analyzing $\cL_{2}$, it is useful to note that when 
$r$ is squarefree and $(a,r)=p$,
\begin{align*}
    \sum_{\substack{ N < n \le 2N \\ n\equiv a \pmod r }}
     \varpi*\varpi(n) 
     = & 2\log p 
     \sum_{\substack{\frac{N}{p} <m\le \frac{2N}{p} \\
		      m\equiv \frac{a}{p} \pmod{\frac{r}{p}} \\
		      }} 
	   \varpi(m) \\
     = & \frac{2N}{\phi(r)} \frac{(\log p)\phi(p)}{p} + 
         O(E^*(N/p,r/p)).	  
\end{align*}
    
When $r$ is squarefree and $p$ is a prime dividing $r$,
we define
\begin{equation} \label{E:OmegadpDef}
    \Omega^*_{r,p}(\cH) = 
       \{a \in \bZ_r: (a,r)=p \text{ and }
	    P(a;\cH)\equiv 0 \pmod r\}.
\end{equation}
Let $\nu^*_{r,p}=\nu^*_{r,p}(\cH)$ be the cardinality of $\Omega^*_{r,p}(\cH)$. 

We take $d_1=p,d_2=r/p$ in \eqref{E:CRT}, and we see that
$\Omega_{r,p}^*(\cH)$ is the image of the set 
$\{0\}\times \Omega_{r/p}^*$ under the isomorphism $\xi$ of 
\eqref{E:CRT}. Therefore
\begin{equation*}
\nu^*_{r,p} = \nu^*_{r/p}.
\end{equation*} 

Using the above information, we find that
\begin{align} \label{E:L22Est1}
    \cL_{2} & =
     \sum_{d,e} \lambda_{d,\ell_1} \lambda_{e,\ell_2} 
        \sum_{p|[d,e]} \sum_{a\in \Omega^*_{[d,e],p}}
	  \sum_{\substack{ N < n \le 2N \\ n\equiv a \pmod {[d,e]} }}
             \varpi*\varpi(n) \\
        & =
           2N \sum_{d,e}  \frac{\lambda_{d,\ell_1}\lambda_{e,\ell_2}}{\phi([d,e])}
		  \sum_{p|[d,e]} \nu^*_{[d,e]/p} \frac{(\log p) \phi(p)}{p}
		  + O\left( (\log N)^{4k} \cE_2 \right),
		  \notag
\end{align}
where
\begin{equation*}
    \cE_2=  \sumflat_{r< R^2} 3^{\omega(r)} 
	      \sum_{\substack{p|r \\ p<R}} 
		     \nu^*_{r/p} E^*(N/p,r/p) (\log p).
\end{equation*}
Upon writing $r=pm$ and changing the order of summation,
we find that
\begin{equation*}
    \cE_2 \le 3 \sum_{p<R} \log p
               \sumflat_{m<R^2/p} 3^{\omega(m)}\nu^*_m E^*(N/p,m).
\end{equation*}
By Lemma \ref{L:HR35}, the inner sum is
$\ll (N/p) (\log N/p)^{-4k-2} \ll (N/p) (\log N)^{-4k-2}$.
Summing over $p$, we get 
\begin{equation*}
   \cE_2 \ll N(\log N)^{-4k-1}.
\end{equation*}
Therefore
\begin{equation} \label{E:L22Est2}
 \cL_{2} = 2N \sum_{d,e}  
               \frac{\lambda_{d,\ell_1}\lambda_{e,\ell_2}}{\phi([d,e])}
		  \sum_{p|[d,e]}  \frac{ \nu^*_{[d,e]/p}(\log p)\phi(p)}{p} 
		  + O(N/\log N).
\end{equation}

Now we turn our attention to $\cL_1$. 
From our definitions and \eqref{E:R1}, we have
\begin{align} \label{E:L1Est1}
\cL_1 = & \sum_{d,e} \lambda_{d,\ell_1} \lambda_{e,\ell_2} 
	     \sum_{a\in \Omega^*_{[d,e]}} 
	     \sum_{\substack{ N \le n < 2N \\ n\equiv a \pmod{[d,e]} }}
	        \varpi*\varpi(n)    \\
    = & N \sum_{d,e} 
	   \frac{\lambda_{d,\ell_1} \lambda_{e,\ell_2} \nu^*_{[d,e]}}
	            {\phi([d,e])}
	   \left(\log N+ C_0 -2  \sum_{p|[d,e]} \frac{\log p}{p} \right)
       + O(\cE_1), \notag
\end{align}
where  
\begin{equation*}
   \cE_1= (\log R)^{4k} \sumflat_{r<R^2} 3^{\omega(r)} \nu^*_r E_2^*(N,r).
\end{equation*}
By Lemma \ref{L:HR35}, $\cE_1 \ll N/\log N$.

Combining our estimates for $\cL_1$ and $\cL_{2}$, we find that
\begin{equation} \label{E:LEst-S5}
\cL=N (\log N+C_0) S_1 -2NS_2 + 2N S_3 + O(N/\log N),
\end{equation}
where
\begin{align*}
S_1= & \sum_{d,e} \frac{\lambda_{d,\ell_1}\lambda_{e,\ell_2}\nu^*_{[d,e]}}
	    {\phi([d,e])} , \\
S_2= & \sum_{d,e} \frac{\lambda_{d,\ell_1} \lambda_{e,\ell_2}\nu^*_{[d,e]}}
	    {\phi([d,e])}
	 \sum_{p|[d,e]} \frac{\log p}{p},\\
S_3= & \sum_{d,e} \frac{\lambda_{d,\ell_1} \lambda_{e,\ell_2}}{\phi([d,e])}
	  \sum_{p|[d,e]} \nu^*_{[d,e]/p} \frac{(\log p)\phi(p)}{p}.\\
\end{align*}

We have already encountered the sum $S_1$; it is the same as the
sum $S$ defined in \eqref{E:SDef-Thm6}. 
From  \eqref{E:SEst}, we
see that
\begin{equation} \label{E:S1EstFinal}
    S_1= 
    \binomial{\ell_1+\ell_2+2}{\ell_1+1} \fS(\cH)
	  \frac{(\log R)^{k+\ell_1+\ell_2+1}}{(k+\ell_1+\ell_2+1)!} 
	  + O\left(
	       \beta(\cH)\fS(\cH)^2 (\log R)^{k+\ell_1+\ell_2}
		 \right).
\end{equation}
 
Of the remaining two sums, $S_3$ is more important, so we concentrate 
on it first. We begin by interchanging the order of summation in $S_3$;
this yields 
\begin{equation} \label{E:S3Decomp1}
    S_3= \sum_p \frac{\log p}{p} U(p),
\end{equation}
where
\begin{equation} \label{E:UDef}
    U(p)= \sum_{\substack{d,e\\p|[d,e]}} 
	  \frac{\lambda_{d,\ell_1}\lambda_{e,\ell_2} \nu^*_{[d,e]/p}}
		 {\phi([d,e]/p)}.
\end{equation}

We decompose $U(p)$ as 
\begin{equation} \label{E:UDecomp}
    U(p)=U_1(p)+U_2(p)+U_3(p),
\end{equation}
where
\begin{align*}
    U_1(p,\ell_1,\ell_2) & = \sum_{\substack{d,e\\p|d, p\nmid e}} 
    \frac{\lambda_{d,\ell_1}\lambda_{e,\ell_2} \nu^*_{[d,e]/p}}
		    {\phi([d,e]/p)},\\
    U_2(p,\ell_1,\ell_2) & = \sum_{\substack{d,e\\p\nmid d, p|e}} 
	\frac{\lambda_{d,\ell_1}\lambda_{e,\ell_2} \nu^*_{[d,e]/p}}
	  {\phi([d,e]/p)},\\
    U_3(p,\ell_1,\ell_2) & = \sum_{\substack{d,e\\p|d, p|e}} 
	    \frac{\lambda_{d,\ell_1}\lambda_{e,\ell_2} \nu^*_{[d,e]/p}}
		{\phi([d,e]/p)}.\\	  
\end{align*}
Going back to 
\eqref{E:S3Decomp1}, 
we will write
\begin{equation} \label{E:S3Decomp2}
    S_3= S_{3,1}+S_{3,2}+S_{3,3},
\end{equation}
where
\begin{equation*}
    S_{3,i} = \sum_p \frac{\log p}{p} U_i(p,\ell_1,\ell_2).
\end{equation*}
We will ultimately see that each  $S_{3,i}$ corresponds to one of the 
terms in the quantity $T(k,\ell_1,\ell_2)$ defined in the statement of 
Theorem \ref{T:Thm7}. 
More precisely, we will show that when $1\le i \le 3$, 
\begin{equation*}
    S_{3,i} = 
      T_i   \fS(\cH) 
        \frac{(\log R)^{k+\ell_1+\ell_2+2}} {(k+\ell_1+\ell_2+2)!}
     \SE,
\end{equation*}
where
\begin{equation*} 
   T_1= -\binom{\ell_1+\ell_2+3}{\ell_2+1}, \quad
   T_2= - \binom{\ell_1+\ell_2+3}{\ell_1+1}, \quad
   T_3= \binom{\ell_1+\ell_2+2}{\ell_1+1}. 
\end{equation*}

We note that $U_2(p,\ell_1,\ell_2)$ is the same as $U_1(p,\ell_1,\ell_2)$ 
except that the roles of 
$\ell_1,\ell_2$ have been reversed; i.e., 
$U_2(p,\ell_1,\ell_2)=U_1(p,\ell_2,\ell_1)$. 
Accordingly, we will concentrate 
on evaluating $U_1(p,\ell_1,\ell_2)$ and $U_3(p,\ell_1,\ell_2)$. 
For brevity, we will usually write these as $U_1(p)$ and $U_3(p)$.

The evaluations of $U_1(p)$ and $U_3(p)$ will require use 
of the quantity $y^*_{r,\ell}$ defined in 
\eqref{E:y*Def}, 
as well as a new quantity $z^*_{r,p,\ell}$. 
The latter is defined as 
\begin{equation} \label{E:z*Def}
    z^*_{r,p,\ell} = 
    \begin{cases}
    \mu(pr) f^*_1(r) 
     \displaystyle{
             \sumprime_d \frac{\lambda_{drp,\ell}}{f^*(dr)}
	         }
      & \text{ if $r< R/p$ and $(r,A)=1$,}\\
      0 & \text{ otherwise.}
     \end{cases}
\end{equation}
As in Section \ref{S:Thm6}, we use $\sumprime_{\null}$ to denote that 
the sum is over values of the indices that are relatively prime to $A$.
Note that $z^*_{r,p,\ell}=0$ if $(p,r)\ne 1$. 
On the other hand, the condition $p|A$ (i.e., $\nu^*_p=0$) does not imply
that $z^*_{r,p,\ell}=0$. However, one can easily 
show that if $p\nmid A$, then
\begin{equation} \label{E:zyRelation}
    z^*_{r,p,\ell} = \left( \frac{p-1}{p-\nu_p} \right) y^*_{rp,\ell}.
\end{equation}

We now give three lemmas that we will use for the 
evaluation of $S_1$ and $S_3$.

\begin{lemma} \label{L:Uyz} If $p<R$, then
    \begin{align}
	U_1(p) & =  - \sumprime_{\substack{r\\(r,p)=1}} 
	                 \frac{z^*_{r,p,\ell_1}y^*_{r,\ell_2}}{f_1^*(r)}
		     - \frac{\nu^*_p}{p-1}
		      \sumprime_{\substack{r\\(r,p)=1}} 
			 \frac{z^*_{r,p,\ell_1} z^*_{r,p,\ell_2}}{f_1^*(r)},
			 \text{ and }
	 \label{E:U1yz}\\
	 U_3(p) & =  \sumprime_{\substack{r\\(r,p)=1}} 
			\frac{z^*_{r,p,\ell_1} z^*_{r,p,\ell_2}}{f_1^*(r)}
	  \label{E:U3yz}.
     \end{align}
\end{lemma}

\begin{proof} The sum $U_1(p)$ may be written as 
\begin{align} 
    U_1(p) & = \sumprime_{\substack{ d,e \\ p\nmid e}} 
                \frac{\lambda_{dp,\ell_1}\lambda_{e,\ell_2}}{\phi([d,e])} 
                \nu^*_{[d,e]} 
	     = \sumprime_{\substack{ d,e \\ p\nmid e}} 
			       \frac{\lambda_{dp,\ell_1}\lambda_{e,\ell_2}}{f^*([d,e])} 
	     = \sumprime_{\substack{ d,e \\ p\nmid e}} 
                \frac{\lambda_{dp,\ell_1}\lambda_{e,\ell_2}}{f^*(d)f^*(e)} 
		 \sum_{\substack{ r|d \\ r|e}} f_1^*(r) 
		 \label{E:U1Form1} \\
	   & = \sumprime_{\substack{r \\ (r,p)=1}} f_1^*(r)
	        \left(\sum_{d} 
		        \frac{\lambda_{drp,\ell_1}}{f^*(dr)} 
			   \right)
		\left(\sum_{\substack{e \\ p\nmid e}} 
			\frac{\lambda_{er,\ell_2}}{f^*(er)} 
			     \right).
	    \notag
\end{align}
In the last expression, the first sum in parentheses is 
$\mu(pr)z^*_{r,p,\ell_1}/f_1^*(r)$. The innermost sum is
\begin{equation*}
    \sumprime_{e} \frac{\lambda_{er,\ell_2}}{f^*(er)} 
     - \sumprime_{\substack{e\\p|e}} \frac{\lambda_{er,\ell_2}}{f^*(er)}
	= \frac{\mu(r) y^*_{r,\ell_2}}{f_1^*(r)} 
       - \sumprime_{\substack{e\\p|e}} \frac{\lambda_{er,\ell_2}}{f^*(er)}.
\end{equation*}
We claim that 
\begin{equation} \label{E:Innermost}
    \sumprime_{\substack{e\\p|e}} \frac{\lambda_{er,\ell_2}}{f^*(er)}
    =  \frac{\nu^*_p \mu(pr) z^*_{r,p,\ell_2}}{(p-1)f_1^*(r)}.
\end{equation}
If $\nu^*_p=0$, then both sides of \eqref{E:Innermost} are $0$. 
If $\nu^*_p\ne 0$, then
\begin{equation*}
    \sumprime_{\substack{e\\p|e}} \frac{\lambda_{er,\ell_2}}{f^*(er)}
   =\sumprime_{e} \frac{\lambda_{epr,\ell_2}}{f^*(epr)}
    =\frac{\mu(pr) z^*_{r,p,\ell_2}} {f_1^*(r) f^*(p)},
\end{equation*}
and \eqref{E:Innermost} follows again. 
 
Going back to \eqref{E:U1Form1}, we find that
\begin{equation*} 
U_1(p)= \sumprime_{\substack{r \\ (r,p)=1}} f_1^*(r)
	  \left( \frac{\mu(rp) z^*_{r,p,\ell_1}}{f_1^*(r)} \right)
	  \left( 
	  \frac{\mu(r) y^*_{r,\ell_2}}{f_1^*(r)} -
		  \frac{\mu(rp) z^*_{r,p,\ell_2}\nu^*_p}{f_1^*(r)(p-1)} 
		  \right),
\end{equation*}
and \eqref{E:U1yz} follows.

For $U_3(p)$, observe that
\begin{align*}
    U_3(p) = & \sumprime_{d,e} 
     \frac{\lambda_{dp,\ell_1} \lambda_{ep,\ell_2}}{f^*([d,e])}
       = \sumprime_{d,e} \frac{\lambda_{dp,\ell_1} \lambda_{ep,\ell_2}}{f^*(d) f^*(e)}
	 \sum_{\substack{r|d\\r|e}} f_1^*(r) 
    \\
    = & \sumprime_{\substack{r\\(r,p)=1}} f_1^*(r) 
	  \left( \sumprime_d \frac{\lambda_{drp,\ell_1}}{f^*(dr)} 
		    \right)
	   \left( \sumprime_e \frac{\lambda_{erp,\ell_1}}{f^*(er)} 
				       \right)
    \\
    = & \sumprime_{\substack{r\\(r,p)=1}} \frac{z^*_{r,p,\ell_1} z^*_{r,p,\ell_2}} {f_1^*(r)},
\end{align*}
and this yields \eqref{E:U3yz}.
\end{proof}

\begin{lemma}  \label{L:z*r}
If $r< R/p$ and $(r,A)=1$, then
    \begin{align}  
	z^*_{r,p,\ell} = 
	\mu^2(rp) &
	 \frac{\fS(\cH)}{(\ell+1)!} 
	  \left( \frac{p-1}{p-\nu_p} \right)
	  (\log R/rp)^{\ell+1} \label{E:z*r}\\
	  & + O\left(\mu^2(rp)
	      \frac{\rho(rp)rp}{\phi(rp)} \fS(\cH)
		   (\log 2R/rp)^{\ell}
		 \right). \notag
	\end{align}
\end{lemma} 

We remark that the error term could be simplified;
it is obvious that 
$$\frac{\rho(rp)rp}{\phi(rp)} \ll \frac{\rho(r)r}{\phi(r)}.$$ 
However, we prefer to write it as above to emphasize the connection
between $y^*_{r,\ell}$ and $z^*_{r,p,\ell}$. In fact,
this lemma follows immediately from \eqref{E:zyRelation} and 
\eqref{E:y*2} when $\nu^*_p \ne 0$.  However, the following
argument works whether or not $\nu^*_p=0$.

\begin{proof} 
The result is trivial is $rp$ is not squarefree,
becuase both sides of  \eqref{E:z*r} are 0 in this case.
For the rest of this proof, we assume that $rp$ is 
squarefree. Note that this assumption implies that $(r,p)=1$.

We start by observing that
\begin{align*}
   \frac{ \mu(rp) z^*_{r,p,\ell} }{f_1^*(r)} & =
       \sumprime_d \frac{\lambda_{drp,\ell}}{f^*(dr)} 
       = \sumprime_d \frac{\mu(drp)}{f^*(dr)} f(drp) 
	    \sum_t \frac{y_{drpt,\ell}}{f_1(drpt)}
	\\
       & = \frac{\mu(rp)f(rp)}{f^*(r)f_1(rp)} 
	    \sumprime_d \frac{\mu(d) f(d)}{f^*(d)} 
	     \sum_t \frac{y_{rpdt,\ell}}{f_1(dt)}
	  \\
	& = \frac{\mu(rp)f(rp)}{f^*(r)f_1(rp)}
	      \sum_{\substack{m\\(m,rp)=1}} 
		\frac{y_{rpm,\ell}}{f_1(m)}
		  \sumprime_{d|m} \frac{\mu(d)f(d)}{f^*(d)} 
	   \\
	  & = \frac{\mu(rp)f(rp)}{f^*(r)f_1(rp)}
	      \sum_{\substack{m\\(m,rp)=1}} 
		\frac{y_{rpm,\ell}}{\phi(m)}.
\end{align*}
In the last line, we have used the relation \eqref{E:y*-Formula1}.
If we also use \eqref{E:y*-Formula2}, we find that  
\begin{align} \label{E:z*1}
    z^*_{r,p,\ell} & = 
	 \frac{f(r)f_1^*(r) f(p)}{f_1(r)f^*(r)f_1(p)} 
	   \sum_{\substack{m< R/rp \\ (m,rp)=1}} 
	      \frac{y_{rpm,\ell}}{\phi(m)} \\
	  & = \frac{\fS(\cH)}{\ell!} 
		\frac{rp}{\phi(rp)} 
		  \frac{(p-1)}{(p-\nu_p)} 
		  \sum_{\substack{m< R/rp \\ (m,rp)=1}} 
			       \frac{\mu^2(m)}{\phi(m)}
				 (\log R/rpm)^{\ell}.
		\notag
\end{align}
We then use \eqref{E:PhiSum} to complete the proof.
\end{proof}

\begin{lemma} \label{L:PrimeSum}
If $a,b$ are non-negative integers, then
\begin{align*}
    \sum_{p< R} \frac{(\log p)^{a+1} (\log R/p)^b}{p}
    = & \frac{a!b!}{(a+b+1)!} (\log R)^{a+b+1}
       + O_{a,b}((\log R)^{a+b}) , \text{ and}\\
      \sum_{p< R} \frac{ (\log p)^{a+1} (\log R/p)^b }{p^2}
       & \ll_a (\log R)^{b}.
\end{align*}
\end{lemma}

\begin{proof} Let $E(u)$ be defined by the relation
    \begin{equation*}
       \sum_{p\le u} \frac{\log p}{p}  = \log u + E(u).
     \end{equation*}
It is well-known that $E(u) \ll 1$.  The first sum in the lemma is
\begin{align*}
  \sum_{p< R} &  \frac{(\log p)^{a+1} (\log R/p)^b}{p}= \\
   & \int_{1}^R (\log u)^{a} (\log R/u)^b \frac{du}{u} +  
	 \int_1^R (\log u)^{a} (\log R/u)^b dE(u).
\end{align*}
By Lemma \ref{L:BetaIntegral}, the first integral is 
\begin{equation*}
    \frac{a!b!}{(a+b+1)!} (\log R)^{a+b+1}. 
\end{equation*}
Using integration by parts, we see that the second integral is
\begin{equation*}
    \int_1^R E(u) 
      \frac{d}{du} \left\{(\log u)^a (\log R/u)^b\right\}  \, du
      \ll_{a,b} (\log R)^{a+b}.
\end{equation*}

This proves the first statement. The second statement is easier; we 
simply note that
\begin{equation*}
    \sum_{p< R} \frac{  (\log p)^{a+1} (\log R/p)^b }{p^2}
    \ll (\log R)^b \sum_{p} \frac{ (\log p)^{a+1} }{p^2} \ll_a (\log R)^b.
\end{equation*}
\end{proof}

{\it Evaluation of $S_{3,3}$}.
From Lemmas \ref{L:Uyz} and \ref{L:z*r}, we see that
\begin{equation} \label{E:U3Est1}
    U_3(p)=
     \frac{\fS(\cH)^2}{(\ell_1+1)!(\ell_2+1)!} 
     \left( \frac{p-1}{p-\nu_p} \right)^2
      V_3(p)
      + 
      O\left( \fS(\cH)^2 (\log R)^{\ell_1+\ell_2+1} W(p) \right),
\end{equation}
where
\begin{align}
     V_3(p) & = \sumprime_{\substack{r< R/p\\ (r,p)=1}}
               \frac{\mu^2(r)}{f_1^*(r)} (\log 
               R/rp)^{\ell_1+\ell_2+2}, \text{ and}
	       \label{E:V3pDef}\\
     W(p) & = \sumprime_{r<R/p} \frac{\mu^2(r)\rho(r) r}{f_1^*(r)\phi(r)}.
                \label{E:WpDef}
\end{align}

$W(p)$ is majorized by the sum $W$ defined in \eqref{E:WDef}, and,
using \eqref{E:WEstFinal}, we see that 
\begin{equation} \label{E:WpEst}
    W(p) \ll   (\log R)^{k-1}.
\end{equation}
From Lemma \ref{L:ExtendedWirsing-f*}, we see that 
\begin{equation} \label{E:V3pEst}
    V_3(p) = 
      \frac{(\ell_1+\ell_2+2)!}{\fS(\cH)}
      \left( \frac{p-\nu_p}{p-1} \right)
      \frac{ (\log R/p)^{k+\ell_1+\ell_2+1}}{(k+\ell_1+\ell_2+1)!}
      + 
      O\left( \beta(\cH) (\log R)^{k+\ell_1+\ell_2} \right).
\end{equation}
We combine the above estimates for $V_3(p)$ and $W(p)$ with 
\eqref{E:U3Est1} to get 
\begin{align}  
    U_3(p)= 
       \binom{\ell_1+\ell_2+2}{\ell_1+1} \fS(\cH) &
	\left( \frac{p-1}{p-\nu_p} \right)
	\frac{ (\log R/p)^{k+\ell_1+\ell_2+1}}{(k+\ell_1+\ell_2+1)!} 
	 \label{E:U3pEstFinal} \\
	     & + 
	  O\left( \beta(\cH)\fS(\cH)^2 (\log R)^{k+\ell_1+\ell_2} \right).
	   \notag 
\end{align}

We can now finish our estimation of $S_{3,3}$. From our definition
and from \eqref{E:U3pEstFinal}, we get
\begin{align*}
    S_{3,3} = & \sum_{p<R} \frac{\log p}{p} U_3(p)  \\
            = & \binom{\ell_1+\ell_2+2}{\ell_1+1} 
	        \frac{\fS(\cH)}{(k+\ell_1+\ell_2+1)!}
		\sum_{p< R} 
		   \left( \frac{\log p}{p}\right)
		   \left( \frac{p-1}{p-\nu_p} \right)
		    (\log R/p)^{k+\ell_1+\ell_2+1}\\
	      &\phantom{1234567890} 
	        +O\left(
		     \beta(\cH)\fS(\cH)^2 (\log R)^{k+\ell_1+\ell_2}
		     \sum_{p<R} \frac{\log p}{p}  
		     \right).
\end{align*}
Now $(p-1)/(p-\nu_p)= 1 +O(1/p)$, so we may use Lemma \ref{L:PrimeSum} to get 
\begin{align}
    S_{3,3}= & \binom{\ell_1+\ell_2+2}{\ell_1+1} \fS(\cH)
	        \frac{(\log R)^{k+\ell_1+\ell_2+2}}{(k+\ell_1+\ell_2+2)!} 
		\label{E:S33EstFinal} \\
	      & \phantom{01234567890}
	      +O\left( \beta(\cH)\fS(\cH)^2 
	        (\log  R)^{k+\ell_1+\ell_2+1}\right).
		\notag
\end{align}

{\it Evaluation of $S_{3,1}$}. The evaluation of $S_{3,1}$ proceeds 
similarly to the evaluation of $S_{3,3}$, but it is somewhat more
involved. We start by defining
\begin{equation} \label{E:U4Def}
    U_4(p) = \sumprime_{\substack{r\\(r,p)=1}} 
	                 \frac{y^*_{r,\ell_2} z^*_{r,p,\ell_1}}{f_1^*(r)},
\end{equation}
and 
\begin{equation} \label{E:S4Def}
    S_4= \sum_{p<R} \frac{\log p}{p} U_4(p).
\end{equation}
Then \eqref{E:U1yz} may be rewritten as 
\begin{equation} \label{E:U1Est1}
    U_1(p)= - U_4(p) -  \frac{\nu_p^*}{p-1}  U_3(p),
\end{equation}
and we may also write
\begin{equation} \label{E:S31Est1}
    S_{3,1}= - S_4 - \sum_{p<R} \frac{(\log p) \nu^*_p }{p(p-1)} U_3(p).
\end{equation}
From \eqref{E:U3pEstFinal}, we see that
\begin{align} 
    \sum_{p<R} \frac{(\log p)\nu^*_p}{p(p-1)} U_3(p) 
     & \ll
    \beta(\cH)\fS(\cH)^2 (\log R)^{k+\ell_1+\ell_2+1} 
      \sum_p \frac{(\log p) \nu^*_p }{p^2} 
      \label{E:S31Error} \\
      & \ll
      \beta(\cH)\fS(\cH)^2 (\log R)^{k+\ell_1+\ell_2+1}.
      \notag
\end{align}

Now we concentrate on $U_4(p)$ and $S_4$. From 
\eqref{E:y*2} and Lemma \ref{L:z*r}, we see that
\begin{equation} \label{E:U4Est1}
    U_4(p)= 
    \frac{\fS(\cH)^2}{(\ell_1+1)!(\ell_2+1)!} 
	 \left( \frac{p-1}{p-\nu_p} \right)
	  V_4(p)
	  + 
	  O\left( \fS(\cH)^2 (\log R)^{\ell_1+\ell_2+1} W(p) \right),
\end{equation}
where $W(p)$ was defined in \eqref{E:WpDef} and
\begin{equation}\label{E:V4pDef}
    V_4(p)= \sumprime_{\substack{r< R/p\\ (r,p)=1}}
               \frac{\mu^2(r)}{f_1^*(r)} 
	        (\log R/r)^{\ell_2+1}
	        (\log R/rp)^{\ell_1+1}.
\end{equation}

We write $\log R/r = \log p + \log R/rp$ and use the 
binomial theorem to get
\begin{equation*}  
    V_4(p) = 
       \sum_{j=0}^{\ell_2+1} \binomial{\ell_2+1}{j} (\log p)^j
	  \sumprime_{\substack { r < R/p \\ (r,p)=1}}
	    \frac{\mu^2(r)}{f_1^*(r)} (\log R/rp)^{\ell_1+\ell_2+2-j}.
\end{equation*}
We apply Lemma \ref{L:ExtendedWirsing-f*} to the inner sum, and  we get
\begin{align} \label{E:V4Est1}
& V_4(p)    = \\
& \frac{1}{\fS(\cH)}  \left( \frac{p-\nu_p}{p-1} \right)   
	      \sum_{j=0}^{\ell_2+1} \binomial{\ell_2+1}{j}  
	      \frac{(\ell_1+\ell_2+2-j)!}{(k+\ell_1+\ell_2+1-j)!}
		 (\log p)^j (\log R/p)^{k+\ell_1+\ell_2+1-j} \notag \\
	   & \phantom{01234567890}
	     + O(\beta(\cH) (\log 2R)^{k+\ell_1+\ell_2} ). \notag
\end{align}  
Using this together with \eqref{E:U4Est1} and \eqref{E:WpEst} gives
\begin{equation} \label{E:U4Est2}
    U_4(p) = \frac{\fS(\cH)}{(\ell_1+1)!(\ell_2+1)!}  
	 U_5(p)  
	    + O(\beta(\cH)\fS(\cH)^2 (\log R)^{k+\ell_1+\ell_2}),  
\end{equation}
where 
\begin{equation} \label{E:U5Def}
    U_5(p) =
      \sum_{j=0}^{\ell_2+1} \binomial{\ell_2+1}{j}  
	      \frac{(\ell_1+\ell_2+2-j)!}{(k+\ell_1+\ell_2+1-j)!}
		 (\log p)^j (\log R/p)^{k+\ell_1+\ell_2+1-j}.
\end{equation}

For future reference, we note the crude estimate
\begin{equation} \label{E:U1Crude}
    U_1(p) \ll \beta(\cH)\fS(\cH)^2 (\log R)^{k+\ell_1+\ell_2+1}
\end{equation}
that is implicit in the combination of \eqref{E:U1Est1},
\eqref{E:U4Est2}, \eqref{E:U5Def}, and \eqref{E:U3pEstFinal}.

Using \eqref{E:S4Def} and \eqref{E:U4Est2}, we see that
\begin{equation} \label{E:S4Est1}
    S_4=  \frac{\fS(\cH)}{(\ell_1+1)!(\ell_2+1)!} 
           \sum_{p< R} \frac{\log p}{p}
	 U_5(p) 
	 + O(\beta(\cH)\fS(\cH)^2 (\log R)^{k+\ell_1+\ell_2+1}).
\end{equation}
We apply 
Lemma \ref{L:PrimeSum} to get 
\begin{align*}
     \sum_{p< R} &
	\left(\frac{\log p}{p}\right) 
	  (\log p)^j (\log R/p)^{k+\ell_1+\ell_2+1-j} \\
	= &
	 \frac{j!(k+\ell_1+\ell_2+1-j)!}{(k+\ell_1+\ell_2+2)!}
	    (\log R)^{k+\ell_1+\ell_2+2}   
	     + O((\log R)^{k+\ell_1+\ell_1+1}).
\end{align*}
Using this in \eqref{E:S4Est1} gives
\begin{align} \label{E:S4Est2}
    S_4 =  \fS(\cH)  &
     \frac{(\log R)^{k+\ell_1+\ell_2+2}}{(k+\ell_1+\ell_2+2)!}
      \sum_{j=0}^{\ell_2+1} 
	\binomial{\ell_1+\ell_2+2-j}{\ell_2+1-j}  \\
	& + O(\beta(\cH) \fS(\cH)^2 (\log R)^{k+\ell_1+\ell_2+1} ).\notag
\end{align}

To treat the sum of binomial coefficients in the above, we make a change 
of variables $j=\ell_2+1-i$. The sum then becomes
\begin{equation} \label{E:BI}
    \sum_{i=0}^{\ell_2+1} \binom{\ell_1+1+i}{i}
    = \sum_{i=0}^{\ell_2+1} 
       \left\{ \binom{\ell_1+2+i}{i} - \binom{\ell_1+1+i}{i-1} \right\},
\end{equation}
provided we make the usual convention that
\begin{equation*}
    \binom{\ell_1+1}{-1} =0.
\end{equation*}
The sum on the right-hand side of \eqref{E:BI} is telescoping, so
\begin{equation*}
  \sum_{i=0}^{\ell_2+1} \binom{\ell_1+1+i}{i} =
  \binom{\ell_1+\ell_2+3}{\ell_2+1}.
\end{equation*}
Putting this information into \eqref{E:S4Est2} gives our final 
estimate for $S_{4}$; i.e.,
\begin{align} 
 S_{4} =  \binomial{\ell_1+\ell_2+3}{\ell_2+1} &\fS(\cH)  
     \frac{(\log R)^{k+\ell_1+\ell_2+2}}{(k+\ell_1+\ell_2+2)!}
       \label{E:S4EstFinal} \\
	& + O(\beta(\cH) \fS(\cH)^2 (\log R)^{k+\ell_1+\ell_2+1} ).\notag
\end{align} 

From this, together with \eqref{E:S31Est1} and \eqref{E:S31Error}, we get 
\begin{align} 
 S_{3,1} = -  \binomial{\ell_1+\ell_2+3}{\ell_2+1} & \fS(\cH)
     \frac{(\log R)^{k+\ell_1+\ell_2+2}}{(k+\ell_1+\ell_2+2)!}
     \label{E:S31EstFinal} \\
	& + O(\beta(\cH) \fS(\cH)^2 (\log R)^{k+\ell_1+\ell_2+1} ).\notag
\end{align}

As we noted earlier, $S_{3,2}$ is the same as $S_{3,1}$ with the 
roles of $\ell_1$ and $\ell_2$ reversed. Therefore
\begin{align} 
    S_{3,2} = -  \binomial{\ell_1+\ell_2+3}{\ell_1+1} & \fS(\cH)
	\frac{(\log R)^{k+\ell_1+\ell_2+2}}{(k+\ell_1+\ell_2+2)!}
	\label{E:S32EstFinal} \\
	   & + O(\beta(\cH) \fS(\cH)^2 (\log R)^{k+\ell_1+\ell_2+1} ).\notag
\end{align} 

Combining \eqref{E:S31EstFinal},\eqref{E:S32EstFinal}, and 
\eqref{E:S33EstFinal} gives
\begin{align}
S_{3} = T(k,\ell_1,\ell_2) & \fS(\cH)
    \frac{(\log R)^{k+\ell_1+\ell_2+2}}{(k+\ell_1+\ell_2+2)!}
    \label{E:S3EstFinal} \\
       & + O(\beta(\cH) \fS(\cH)^2 (\log R)^{k+\ell_1+\ell_2+1} ),\notag
\end{align}
where $T(k,\ell_1,\ell_2)$ is as defined in Theorem \ref{T:Thm7}.

Finally, we will quickly dispatch $S_2$. We rewrite this sum as
\begin{equation*} \label{E:S2Decomp}
    S_2= \sumprime_{d,e} 
	 \frac{\lambda_{d,\ell_1} \lambda_{e,\ell_2}}{f^*([d,e])}
	    \sum_{p|[d,e]}\frac{\log p}{p} 
	 =\sumprime_p \frac{\log p}{p f^*(p)} U(p),
\end{equation*}
where $U(p)$ was defined in \eqref{E:UDef}. 
We employ the crude estimate
\begin{equation*}
     U(p) \ll \fS(\cH)^2 \beta(\cH) (\log R)^{k+\ell_1+\ell_2+1}.
\end{equation*}
This is easily seen by combining \eqref{E:UDecomp}, 
\eqref{E:U1Crude}, \eqref{E:U3pEstFinal}, and using
the symmetry between $U_1(p)$ and $U_2(p)$. 
The sum 
\begin{equation*}
    \sumprime_{p\le R} \frac{\log p}{pf^*(p)} 
\end{equation*}
is $\ll 1$. Combining the above gives the bound
\begin{equation} \label{E:S2EstFinal}
    S_2 \ll \fS(\cH)^2 \beta(\cH) (\log R)^{k+\ell_1+\ell_2+1}.
\end{equation}

The proof of Theorem \ref{T:Thm7} is completed by combining
\eqref{E:LEst-S5} together with the final estimates for $S_1,S_2,S_3$,
which are \eqref{E:S1EstFinal}, \eqref{E:S2EstFinal}, and 
\eqref{E:S3EstFinal} 
respectively.

\section{Proofs of Theorems \ref{T:GPY1} through \ref{T:E2GapsEH}} 
\label{S:OtherThms}

Let $\cH=\{h_1 , h_2 ,\ldots , h_k\}$ be an arbitrary admissible 
$k$-tuple. Without loss of generality, we may specify that 
\begin{equation*}
h_1 < h_2 < \ldots < h_k.
\end{equation*}
It is also useful to assume that 
\begin{equation} \label{E:HBound}
 h_k \le \log N.
\end{equation}
With this hypothesis, we see from Lemma \ref{L:HBound} that the error terms in 
Theorems \ref{T:Thm5}, \ref{T:Thm6}, \ref{T:Thm7} satisfy
\begin{equation*}
    \beta(\cH)\fS(\cH)/\log N 
    \ll (\log\log \log N)^{b_k+1}/\log N
    \ll (\log\log N)/\log N.
\end{equation*}

Consider the sum
\begin{equation} \label{E:cS1Def}
\cS_1:=\sum_{N< n \le 2N}  
 \left\{  \sum_{h\in \cH} \varpi(n+h)  \,-\, (\log 3N) \right\} 
    \left( 
	    \sum_{\ell=0}^L b_\ell (\log R)^{-\ell} \Lambda_{R}(n;\cH,\ell) 
      \right)^2.
\end{equation}

For a given $n$, the sum inside the brackets is non-positive 
unless there are at least two distinct values, $h_i,h_j\in \cH$ such
that $n+h_i,n+h_j$ are primes. 
Consequently, if we can show that
the sum in \eqref{E:cS1Def} is 
$\gg N\fS(\cH)(\log R)^{k+1}$, then we can
conclude that $\liminf_{n\to\infty} ( p_{n+1}-p_n )\le h_k -h_1$. 

Expanding the square in \eqref{E:cS1Def}, we see that 
\begin{equation*}
    \cS_1= 
     \sum_{0\le \ell_1,\ell_2 \le L} 
	b_{\ell_1} b_{\ell_2} (\log R)^{-\ell_1-\ell_2} \cM_1(\ell_1,\ell_2),
\end{equation*}
where
\begin{equation*}
    \cM_1({\ell_1,\ell_2}) = 
    \sum_{N \le n < 2N}
     \left\{ \sum_{h\in \cH} \varpi(n+h) - (\log 3N) \right\}
      \Lambda_{R} (n;\cH,\ell_1)\Lambda_{R} (n;\cH,\ell_2).
\end{equation*}
We assume Hypothesis $BV(\theta)$, and we use Theorems \ref{T:Thm5}
and \ref{T:Thm6} with $R=N^{(\theta-\epsilon)/2}$ to get
\begin{align*} 
    \cM_1({\ell_1,\ell_2}) & \asymptotic   
      \binomial{\ell_1+\ell_2+2}{\ell_1+1} N \fS(\cH) 
       k \frac{(\log R)^{k+\ell_1+\ell_2+1}}{(k+\ell_1+\ell_2+1)!}  \\
       & \phantom{1234567890}  -
    \binomial{\ell_1+\ell_2}{\ell_1} N \fS(\cH) 
	\frac{(\log R)^{k+\ell_1+\ell_2} \log N}{(k+\ell_1+\ell_2)!} 
	\\
     & \asymptotic 
	 N\fS(\cH)  (\log R)^{k+\ell_1+\ell_2} (\log N)
	 (m(k,\ell_1,\ell_2,\theta) -\epsilon')
\end{align*}
where 
\begin{equation} \label{E:mDef}
    m(k,\ell_1,\ell_2,\theta)=
     \binomial{\ell_1+\ell_2}{\ell_1} \frac{1}{(k+\ell_1+\ell_2)!}
	\left(
	  \frac{k(\ell_1+\ell_2+1)(\ell_1+\ell_2+2)}
		 {(k+\ell_1+\ell_2+1)(\ell_1+1)(\ell_2+1)}
	   \frac{\theta}{2}
	  -1
	  \right),
\end{equation}
and $\epsilon'=\epsilon'(k,\ell_1,\ell_2,\epsilon)$
goes to 0 as $\epsilon$ goes to $0$.

Define ${\bb}=(b_0,b_1,\ldots,b_L)$. 
Then (we suppress the $\epsilon'$ term)
\begin{align} 
    \cS_1^*(N,\cH,\theta,{\bb}) & :=
    \frac{\cS_1}{N\fS(\cH) (\log R)^k \log N} 
    \label{E:S1*Eq}\\
    & \asymptotic \sum_{0\le \ell_1,\ell_2 \le L} 
	b_{\ell_1} b_{\ell_2} m(k,\ell_1,\ell_2,\theta) \notag\\
    & = \bb^T \bM \bb, \notag
\end{align}
where $\bM=\bM(k,\theta)$ is the matrix
\begin{equation*}
    \bM = \left[ m(k,i,j,\theta) \right]_{0\le i,j \le L}.
\end{equation*}
Our goal is to pick $\bb$ to make $\cS_1^*>0$ for a given $\theta$
and minimal $k$. This is easily determined by picking
$\bb$ to be an eigenvector of the matrix $\bM$ with eigenvalue
$\lambda$, in which case
\begin{equation*}
    \cS_1^*  \asymptotic \bb^T \lambda \bb = 
     \lambda \sum_{i=0}^L b_i^2.
\end{equation*}
This will be positive provided $\lambda$ is positive.
We conclude that $\cS_1^*>0$ if $\bM$ has a positive eigenvalue
and $\bb$ is chosen to be the corresponding eigenvector. 

With $k=6$ and $L=1$, we find that 
\begin{equation*}
     \bM= \frac{1}{8!}
       \left[
	  \begin{matrix} 
	       48\theta-56 & 9\theta-8 \\
	       9\theta-8 & 2\theta-2 \\
	   \end{matrix}
	 \right].
\end{equation*} 
The determinant of $8!\bM$ is $15\theta^2-64\theta +48$, which is 
negative if $4(8-\sqrt{19})/15 < \theta \le 1$. Since the determinant
is the product of the eigenvalues, we conclude that $\bM$ has a 
positive eigenvalue for $\theta$ in this range. 
Consequently, if $\cH$ is an admissible $6$-tuple, then there
are infinitely many $n$ such that at least two of the numbers
$n+h_1, \ldots, n+h_6$ are prime.
We complete the proof of the second part of Theorem 2 by taking
\begin{equation*}
    \cH=\{7,11,13,17,19,23\}.
\end{equation*} 
$\cH$ is admissible because for $p\le 5$, none of the elements in 
$\cH$ are divisible by $p$, and for $p\ge 7$, there are not enough 
elements to cover all of the residue classes mod $p$. 

To prove the first part of Theorem \ref{T:GPY2}, we again use 
\eqref{E:S1*Eq}; however, we use the trivial choice
$b_\ell=1$ for some specific $\ell$,
and $b_i=0$ for all other $i$. Then
\begin{equation*}
    \cS^*_1 \asymptotic m(k,\ell,\ell,\theta) =
	   \binomial{2\ell}{\ell} \frac{1}{(k+2\ell)!}
	\left(
	  \frac{2k(2\ell+1)}{(k+2\ell+1)(\ell+1)}
	   \frac{\theta}{2}
	  -1
	  \right) - \epsilon'.    
\end{equation*}
The above is positive if 
\begin{equation*}
\theta> 
  \left(\frac{1}{2} + \frac{1}{4\ell+2}\right)
  \left( 1+ \frac{2\ell+1}{k} \right).
\end{equation*}
The right-hand side approaches $1/2$ if $\ell,k\to\infty$ with
$\ell=o(k)$.

The above argument just fails when $\theta = 1/2$.
To remedy this, we modify
\eqref{E:cS1Def} by taking $h$ to be a parameter
to be chosen later, with $h\le \log N$.
We then sum
over all admissible size $k$ subsets $\cH$ of $\{1, \ldots, h\}$.
Specifically, we take 
\begin{equation} \label{E:STLambda}
\ctS_1=\sum_
     {\substack{
	 \cH\subseteq \{1,\ldots,h\} \\ |\cH|=k \\ 
	      \cH \text{ admissible}  }}
   \sum_{N<n\le 2N} 
      \left\{  
	   \sum_{1\le h_0 \le h} \varpi(n+h_0) \,-\, (\log 3N) 
      \right\}
	   \Lambda_{R}^2(n;\cH,\ell).
\end{equation}

We apply Theorems \ref{T:Thm5} and \ref{T:Thm6} to the sum
$\ctS_1$ for those terms when $\cH$ and $\cH\union \{h_0\}$
are both admissible. There may be terms with 
$\cH$ admissible but $\cH\union\{h_0\}$ not admissible;
for these terms we apply the trivial bound
\begin{equation*}
    \sum_{N < n \le 2N} 
      \sum_{1\le h_0 \le h} \varpi(n+h_0) \Lambda_{R}(n;\cH,\ell)^2 
      \ge  0. 
\end{equation*}
We find that
\begin{align} \label{E:S1Twiddle}
\ctS_1 \gtrsim & 
  \binom{2\ell+2}{\ell+1}
 \frac{N(\log R)^{k+2\ell+1}}{(k+2\ell+1)!}  
    \sum_{1\le h_0 \le h}
      \sum_{\substack{\cH\subseteq \{1,\ldots, h\} \\ 
		      |\cH|=k, h_0\in  \cH  }} 
		\fS(\cH) \\
  & + 
     \binom{2\ell}{\ell}
      \frac{N(\log  R)^{k+2\ell}}{(k+2\ell)!}
	\sum_{1\le h_0 \le h}
	 \sum_
	     {\substack{\cH \subseteq \{1,\ldots, h\} \\ 
			      |\cH|=k, h_0\notin  \cH }}
		\fS\left(\cH \union \{h_0\}\right) \notag\\
  & - 
      \binom{2\ell}{\ell}
       \frac{N (\log N)(\log R)^{k+2\ell}  }{(k+2\ell)!} 
	  \sum_
	     {\substack{\cH \subseteq \{1,\ldots, h \} \\|\cH|=k }
		} 
	      \fS(\cH). \notag
\end{align}
We have dropped the condition that $\cH$ is admissible in the above 
sums; we may do so because $\fS(\cH)=0$ when $\cH$ is not admissible.

Now we observe that
\begin{equation*}
    \sum_{1\le h_0 \le h}
	 \sum_{\substack{\cH\subseteq \{1,\ldots, h\} \\ 
			 |\cH|=k, h_0\in  \cH }} 
		   \fS(\cH) 
= k    \sum_{\substack{\cH\subseteq \{1,\ldots, h\} \\ 
			 |\cH|=k }} \fS(\cH)
\asymptotic \frac{kh^k}{k!}.
\end{equation*}

In the above, equality occurs from noting that every relevant
set $\cH$ occurs $k$ times in the initial sum, and the asymptotic
relation is a theorem of Gallagher \cite{Gallagher}.
We also have that
\begin{equation*}
 \sum_{1\le h_0 \le h}
    \sum_{\substack{\cH\subseteq \{1,\ldots, h\} \\ 
				  |\cH|=k, h_0\notin  \cH }}
		    \fS(\cH \cup \{h_0\})
= (k+1) \sum_{\substack{\cH\subseteq \{1,\ldots, h\} \\ 
				  |\cH|=k+1 }}
		    \fS(\cH)
\asymptotic \frac{h^{k+1}}{k!}.
\end{equation*}

Returning to the evaluation of $\ctS_1$, we find that 
\begin{equation*}
\ctS_1 \gtrsim
\binom{2\ell}{\ell}
       \frac{N (\log N)(\log R)^{k+2\ell} h^k  }{k!(k+2\ell)!} 
\tb_1(k,\ell,h)
\end{equation*} 
where 
\begin{equation*}
\tb_1(k,\ell,h)=
 2\cdot \frac{2\ell+1}{\ell+1} \cdot \frac{k}{k+2\ell+1} \cdot 
     \frac{\log R}{\log N} 
  + \frac{h}{\log N} -1.
\end{equation*}

Unconditionally, we may take $\theta=1/2$, so 
$\log R/\log N=1/4-\epsilon$.  We get two primes in some interval
$(n,n+h], N< n \le 2N$ provided $\tb_1(k,\ell,h)>0$. 
This is equivalent to 
\begin{align*}
\frac{h}{\log N} > &
  1- \frac{2k}{k+2\ell+1}\cdot \frac{2\ell+1}{\ell+1} \cdot
	\left( \frac{1}{4} - \epsilon \right) \\
 = & \frac{k+4\ell^2+6\ell+2+4\epsilon(k+2k\ell)}
	   {2(1+\ell)(1+2\ell+k)}.
\end{align*}
On letting $\ell=[\sqrt{k}]$ and taking $k$ sufficiently large,
we see that this is valid with $h/\log N$ arbitrarily small.
This proves Theorem \ref{T:GPY1}.

For the proofs of Theorem \ref{T:E2Gaps} and
Theorem \ref{T:E2GapsEH}, we note that if 
$N < n \le 2N$ then
\begin{equation*}
\varpi*\varpi(n) \le \frac{(\log 3N)^2}{2}.
\end{equation*}
Accordingly, we consider
\begin{align} \label{E:cS2Def}
    \cS_2:=\sum_{N< n \le 2N} &  
     \left\{  \sum_{h\in \cH} \varpi*\varpi(n+h)  \,-\, 
	     \frac{(\log3N)^2}{2} \right\} \times \\
       & \times\left( 
		\sum_{\ell=0}^L b_\ell (\log R)^{-\ell} \Lambda_{R}(n;\cH,\ell) 
	  \right)^2. \notag
\end{align}
The term $n$ contributes a negative amount unless there are two values
$h_i,h_j\in \cH$ such that $n+h_i,n+h_j$ are products of two primes. 
The values of $n$ for which any $n+h$ is a square of a prime
contribute $\ll N^{1/2} (\log N)^{2k+2}$, and this contribution
may be absorbed into the error terms of our estimates.

We assume Hypotheses $BV(\theta)$ and $BV_2(\theta)$, and we 
argue along the same lines as in the proof of Theorem \ref{T:GPY2}.
When $R=N^{(\theta-\epsilon)/2}$, we obtain
\begin{equation*}
\cS_2 =\sum_{0\le \ell_1,\ell_2 \le L} 
	 b_{\ell_1} b_{\ell_2} (\log R)^{-\ell_1-\ell_2} \cM_2(\ell_1,\ell_2),
\end{equation*}
where
\begin{align*}
\cM_2 & \asymptotic 
\fS(\cH) N (\log N)^2 (\log R)^{k+\ell_1+\ell_2} 
(m_2(k,\ell_1,\ell_2,\theta)-\epsilon'), \\
m_2(k,\ell_1,\ell_2,\theta) & = m_{21}+ m_{22} - m_{23},\\
\end{align*}
\begin{align*}
m_{21} & = \binomial{\ell_1+\ell_2+2}{\ell_1+1} 
	\frac{k}{(k+\ell_1+\ell_2+1)!} \frac{\theta}{2},\\
m_{22} & = 2 \left\{ 
	    \binomial{\ell_1+\ell_2+2}{\ell_1+1} 
	   -  \binomial{\ell_1+\ell_2+3}{\ell_1+1} 
	   -  \binomial{\ell_1+\ell_2+3}{\ell_2+1} 
		\right\}
	    \frac{k}{(k+\ell_1+\ell_2+2)!} \frac{\theta^2}{4},\\
m_{23} & = \frac{1}{2} \binomial{\ell_1+\ell_2}{\ell_1} 
	   \frac{1}{(k+\ell_1+\ell_2)!},
\end{align*}
and $\epsilon'=\epsilon'(k,\ell_1,\ell_2,\epsilon)\to 0$ 
as $\epsilon\to 0$.

Let $\bb$ be as defined before. Then
(suppressing the $\epsilon'$ term)
\begin{align*} 
\cS_2^*(N,\cH,\theta,\bb) & := 
     \frac{\cS_2}{N\fS(\cH)(\log R)^k (\log N)^2} 
 \asymptotic \sum_{0\le \ell_1, \ell_2 \le L} 
	 b_{\ell_1}b_{\ell_2} m_2(k,\ell_1,\ell_2,\theta)  \notag \\
& 
= \bb^T \bM_2 \bb,
\end{align*}
where $\bM_2=\bM_2(k,\theta)$ is the matrix
\begin{equation*}
\bM_2=[m_2(k,i,j,\theta)]_{0\le i,j \le L}.
\end{equation*}

We first prove Theorem \ref{T:E2GapsEH}.
As in the proof of Theorem \ref{T:GPY2}, we wish to show 
that there is some $\bb$ such that $\cS_2^* >0$ for a given
$\theta$ and minimal $k$. Taking $k=3$ and $L=1$, we find that
\begin{equation} \label{E:bM2withk=3}
\bM_2 = \frac{1}{480}
      \left[
	 \begin{matrix} 
	      -24\theta^2+60\theta-40 & -7\theta^2+18\theta-10 \\
	      -7\theta^2+18\theta-10 & -2\theta^2 +6\theta-4 \\
	  \end{matrix}
	\right].
\end{equation}
If we take $b(0)=1,b(1)=4$, then we find that 
\begin{equation*}
\bb^T \bM_2 \bb = -\frac{7\theta^2}{30} + \frac{5\theta}{8} - 
\frac{23}{60}.
\end{equation*}
This is positive whenever
\begin{equation*}
\frac{75-\sqrt{473}}{56} < \theta \le 1.
\end{equation*}
Finally, we note that $\cH=\{5,7,11\}$ is
an admissible $3$-tuple, so this completes
the proof of Theorem \ref{T:E2GapsEH}.

We can also prove Theorem \ref{T:E2GapsEH} with a slightly wider range 
of allowable $\theta$ by taking the determinant of the matrix in 
\eqref{E:bM2withk=3}. A numerical calculation shows that this 
determinant has a zero at $\theta= 0.943635 \ldots$.

For the proof of Theorem \ref{T:E2Gaps}, we take 
$k=8, L=2, \theta=1/2-\epsilon$, 
and we find that
\begin{equation*}
  \bM_2 = \frac{1}{14!}
\left[
	\begin{matrix} 
	     -216216 & 8736 & 3458 \\
	      8736 & -364 & 14 \\
	      3458 & 14 & -36 \\
	 \end{matrix}
       \right],
\end{equation*}
With
\begin{equation*}
b(0)=1,b(1)=16, b(2)=16,
\end{equation*}
we find that 
\begin{equation*}
 14!\bb^T \bM \bb = 78760>0.
\end{equation*}
Now $\cH=\{11,13,17,19,23,29,31,37\}$ is an admissible
$8$-tuple, so this completes the proof of Theorem \ref{T:E2Gaps}.

We make one final comment regarding the proofs that make use of 
bilinear forms in $\bb$. 
By taking 
\begin{equation*}
    \sum_{\ell=0}^L b_\ell (\log R)^{-\ell} \Lambda_R(n;\cH,\ell)
\end{equation*}
in the definitions of $\cS_1$ and $\cS_2$, we are in essence
using 
\begin{equation*}
    y_r = 
     \fS(\cH) 
      \sum_{\ell=0}^L 
         \frac{b_\ell}{\ell!}
	  \left(\frac{\log R/r}{\log R}\right)^\ell.
\end{equation*}
In other words, we have essentially replaced $(\log R/r)^\ell$
in \eqref{E:yrChoice} by a polynomial in $\log R/r$.

{\it Acknowledgements:} We thank Tsz-Ho Chan and Yoichi Motohashi
for their comments on this paper.
Part of the work for this paper was done at the 
American Institute of Mathematics,  where Graham was 
visiting in Fall 2004. He thanks them for their hospitality and 
excellent working environment.

\end{document}